%% file: K3CarpetEquations.tex
\documentclass[twoside,12pt, leqno]{amsart}
\usepackage{amsmath,amscd,amsthm,amssymb,amsxtra,latexsym,epsfig,epic,graphics}
\usepackage[matrix,arrow,curve]{xy}
\usepackage{graphicx}
\usepackage{diagrams}
\usepackage{tikz}  
\usepackage{color} 
\usetikzlibrary{arrows,calc} 
\usetikzlibrary{decorations.pathmorphing}
\usetikzlibrary{decorations.markings} 
\usetikzlibrary{decorations.pathreplacing} 
\usetikzlibrary{plothandlers}

\input preamble.tex

\def\AA{{\mathbb A}}

\def\PP{{\mathbb P}}
\def\P{{\mathbb P}}

\def\FF{{\mathbb F}}

\def\name{{\rm name}}

\def\CO{{\mathcal O}}

\newarrow{Iso} -----

\def\init{{\rm in}}

\def\lbracket{{[\kern-1.5pt[}}
\def\rbracket{{]\kern-1.5pt]}}

\makeatletter
\def\Ddots{\mathinner{\mkern1mu\raise\p@
\vbox{\kern7\p@\hbox{.}}\mkern2mu
\raise4\p@\hbox{.}\mkern2mu\raise7\p@\hbox{.}\mkern1mu}}
\makeatother

\usepackage{times}
\newdimen\x \x=12pt

\usepackage{color}

\usepackage[breaklinks,bookmarksopen,bookmarksnumbered,urlcolor=blue]{hyperref}
\hypersetup{colorlinks=true,backref=true,citecolor=blue}

\author[David Eisenbud]{David Eisenbud}
\address{Department of Mathematics, University of California at Berkeley and the Mathematical
Sciences Research Institute, Berkeley, CA 94720, USA}
\email{de@msri.org}

\author[Frank-Olaf Schreyer]{Frank-Olaf Schreyer}
\address{Fachbereich Mathematik, Universit\"at des Saarlandes, Campus E2 4, D-66123 Saar\-br\"ucken, Germany}
\email{schreyer@math.uni-sb.de}

\title{Equations and Syzygies of\\ K3 Carpets and Unions of Scrolls}
\begin{document}

\begin{abstract}
We describe the equations and Gr\"obner bases of some degenerate K3 surfaces associated to rational normal scrolls. These K3 surfaces are members of a class of interesting singular projective varieties we call correspondence scrolls. The ideals of these surfaces are nested in a simple way that allows us to analyze them inductively. We describe explicit Gr\"obner bases and syzygies for these objects over the integers and this lets us treat them in all characteristics simultaneously.  
\end{abstract}

\maketitle

\section*{Introduction} 
\let\thefootnote\relax\footnote{
\noindent AMS Subject Classification:
Primary: 14H99,
Secondary: 13D02, 14H51 \smallbreak
Keywords: K3 Surfaces, Green's Conjecture in positive characteristic, canonical curves, canonical ribbons, K3 carpets.\smallbreak
The first author is grateful to the
National Science Foundation for partial support. This work is a contribution to Project I.6 of the second author within the SFB-TRR 195 "Symbolic Tools in Mathematics and their Application" of the German Research Foundation (DFG).}

Let $S(a,b)$ be the rational normal surface scroll of degree $a+b$ in $\PP^{a+b+1}$ over an arbitrary field $\FF$, that is, the embedding of the projectivised vector bundle $ \PP(\CO_{\PP^{1}}(a)\oplus \CO_{\PP^{1}}(b))$ by the line
bundle $\CO(1)$ (see \cite{EH} for an exposition). A striking theorem of Gallego and Purnaprajna (\cite[Theorem 1.3]{GP}) asserts that there is
a unique \emph{K3 Carpet} that is a double structure on $S(a,b)$; that is, a unique scheme $X(a,b)\subset \PP^{a+b+1}$ whose reduced scheme $X(a,b)_{\rm red}$ is  $S(a,b)$ such that $X(a,b)$ has
degree $2(a+b)$ with $\HH^{1}(\CO_{X(a,b)}) = 0$ and $\omega_{X(a,b)} \cong \CO_{X(a,b)}$ (or, equivalently, with
homogeneous coordinate ring Gorenstein of $a$-invariant 0.) Gallego and Purnaprajna  prove that $X(a,b)$ can be written as a limit of smooth K3 surfaces whose general hyperplane sections are canonical curves of genus $a+b-1$ and gonality $\min(a,b)+2$.

A quick  description of the homogeneous ideal of $X(a,b)$ is  that, for $a,b\geq 2$, it is is generated by  the rank 3
quadrics in the ideal of $S(a,b)$ (Theorem \ref{rank 3}). The goal of this paper is to elucidate the generators of  this ideal,
and those of certain related varieties, in a much more explicit way, similar to the well-known description
of the ideal of $S(a,b)$ as an ideal of $2\times 2$ minors. This enables us to compute explicit Gr\"obner bases
and even resolutions over the integers. 

One of our motivations has to do with Green's conjecture relating the Clifford index of a smooth projective
curve to the length of the linear strand of its free resolution. Deopurkar~\cite{D} has recently proven that all canonical ribbons satisfy Green's conjecture. Since every canonical ribbon of genus $g$ and Clifford index $c$ is the hyperplane section of the K3 carpet $X(c, g-1-c)$ (\cite[Section 8]{BE}), this implies that all K3 carpets satisfy the analogue of Green's conjecture.  One can also hope that K3 carpets could shed some light on the questions of the stability of syzygies raised in~{\cite{DFS} }.

Deopurkar's argument relies on Voisin's theorem \cite{V05} that  canonical curves  lying on sufficiently general K3 surfaces satisfy Green's conjecture. In very recent work, Aprodu, Farkas, Papadima, Raicu and Weyman \cite{AFPRW} have given a far simpler proof of Voisin's theorem based on  the degeneration of K3 surfaces to tangent developable surfaces of rational normal curves. 

It seems natural to hope that there might  also be a  proof based on  K3 carpets, and this would have the advantage that it would automatically treat curves of every Clifford index: indeed, the analogue of Green's Conjecture for $X(a,a)$ (which corresponds to Green's conjecture for general curves) directly implies Green's conjecture for all $X(a,b)$ with $b \le a$, and thus for some curves
of each Clifford index. This is because  a Gr\"obner basis for the ideal of  each $X(a,b)$ with $b < a$ is a subset of that of $X(a,a)$. 

Green's Conjecture is known to fail in some finite characteristics (\cite{Bopp}, \cite{BS18}). Because the Gr\"obner bases we construct are valid over the integers, we are able to tabulate the characteristics of the fields over which the conjecture  fails for K3 carpets of sectional genus up to 15 and thus for canonical ribbons of these genera. The data lead us to conjecture:

\begin{conjecture}\label{Conj0}
 Green's conjecture is true for general curves of genus $g$ over fields of characteristic $p>0$ whenever $p\geq (g-1)/2$.
\end{conjecture}

The evidence for this conjecture is presented in more detail in the last section.

\subsection*{Three examples of K3 Carpets}

\noindent{\bf 1)} $S(1,1) \subset X(1,1)$: Any quartic equation in 4 variables defines
a scheme that has the characteristics of a K3 surface. The scroll $S(1,1)$ is  a smooth quadric surface in 
$\PP^{3}$. The unique double structure $X(1,1)$ is defined by the square of the form defining the quadric.
\medbreak

\noindent{\bf 2)} $S(2,1)\subset X(2,1)$: In suitable coordinates $S(2,1)$ is defined by the
$2\times 2$ minors of the matrix 
$$
\begin{pmatrix}
 x_{0}&x_{1}&y_{0}\\
 x_{1}&x_{2}&y_{1}\\
\end{pmatrix}\,.
$$
The carpet $X(2,1)$ supported on this scroll is the  complete intersection
defined by the $2\times2$ minor in the upper left corner, together with the determinant, of the symmetric matrix
$$
\begin{pmatrix}
 x_{0}&x_{1}&y_{0}\\
x_{1}&x_{2}&y_{1}\\
y_{0}&y_{1}&0
\end{pmatrix}\,.
$$
\medbreak

\noindent{\bf 3)} $S(2,2)\subset X(2,2)$: For a more typical example, we take
$S(2,2)$ to be the scroll defined by the $2\times 2$ minors of
$$
\begin{pmatrix}
 x_{0}&x_{1}&y_{0}&y_{1}\\
 x_{1}&x_{2}&y_{1}&y_{2}\\
\end{pmatrix}
$$
then $X(2,2)$ is defined by the complete intersection of the three quadrics
$$\det 
\begin{pmatrix}
 x_{0}&x_{1}\\
 x_{1}&x_{2}\\
\end{pmatrix},
\det 
\begin{pmatrix}
y_{0}&y_{1}\\
y_{1}&y_{2}\\
\end{pmatrix},
\det 
\begin{pmatrix}
x_{0}+y_{0}&x_{1}+y_{1}\\
x_{1}+y_{1}&x_{2}+y_{2}\\
\end{pmatrix}.
$$
We shall see  other useful representations as well.

\subsection*{What's in this paper}
In Section~\ref{general} below we describe a family of projective schemes we call \emph{correspondence scrolls} that includes the rational normal scrolls, and the degenerate K3 surfaces treated in the
rest of this paper. 
In Section~\ref{geometry}, we give an informal description of the family of degenerate K3 surfaces that depend on a pair of automorphisms of $\PP^1$,
and describe their degeneration to a K3 carpet.

Our main results are in Sections \ref{Sequations} and \ref{resolutions}. In Section~\ref{Sequations}, we give various descriptions of the minimal generators of the ideals of the K3 carpets and certain reducible K3 surfaces, and prove
that these generators form a Gr\"obner basis for a suitable term order. 

In Section~\ref{resolutions}, we study a non-minimal free resolutions of these surfaces
that have simple descriptions valid over the ring of integers. Explicit computation then yields information about the characteristics in which Green's conjecture might fail. 

Finally, in Section~\ref{conjectures}, we formulate two Conjectures about the minimal free resolutions of these
surface, and present the data which give the evidence. In particular, we proof Conjecture \ref{Conj0} for curves of genus $g \le 15$.

\section{Correspondence Scrolls}\label{general} 

Consider disjoint projective spaces
$\PP^{a_i} = \PP(V_i)$, for $i = 1,\dots,m$, embedded  in
$$
  \PP^{N} = \PP(\oplus_i V_i),
$$
and a \emph{correspondence}, that is a subscheme
$\Gamma \subset \prod_i\PP^{a_i}$ (or more generally a multi-homogeneous subscheme of 
$\prod_i\AA^{1+a_i}$). The \emph{correspondence scroll}  $S_\Gamma$ defined by
$\Gamma$ may be described set-theoretically as the union of the planes in 
$\PP^N$ spanned by the sets of points $\{p_1,\dots,p_m\}$ with
$(p_1,\dots,p_m)\in \Gamma$. To $S_{\Gamma}$ scheme-theoretically, we first consider
the set of  of planes
of dimension $m-1$ in $\PP^N$ that are spanned by all sets of points
$\{p_1,\dots,p_m\}$ with $p_i \in \PP^a_i\subset \PP^N$. We consider this 
set as a subvariety of the Grassmannian. As such, it is the image of
the product $\prod_i\PP^{a_i}$. We pull back the tautological bundle of $m-1$-planes
on the Grassmannian to
$\Gamma\subset \prod_i\PP^{a_i}$, and we define 
$S_\Gamma$ to be the image in 
$\PP^N$ of this bundle over $\Gamma$.

For example, the ordinary surface scroll $S(a,b)$ is the result of taking 
$$
m =2, \ a_1 =a,  a_2 = b
$$
and taking $\Gamma$ to be the diagonal in $\PP^1\times \PP^1$ embedded in
$\PP^a\times \PP^b$ as the product of the rational normal curves of degrees $a$
and $b$. The  K3 carpet $X(a,b)$ 
described below is obtained by taking $\Gamma$ to be the image of twice the diagonal of $\PP^1\times \PP^1$, and  the other degenerate K3 surfaces we consider correspond to other divisors of type $(2,2)$ on $\PP^1\times \PP^1$. 

It is not hard to describe correspondence scrolls that have the properties of Calabi-Yau varieties of other dimensions, and to give other interesting singular models. This is the subject a paper in preparation by the first author and Allessio Sammartano \cite{EiSa}

In the next section we concentrate on the family of degenerate K3 surfaces.


\section{Degenerate K3 surfaces from rational normal scrolls: Geometry
}\label{geometry}
In this section we sketch the  geometry of the reducible surfaces whose equations we will study. 


Fix positive integers $a,b$, and consider 2-dimensional rational normal scrolls of type $(a,b)$ in $\PP^{a+b+1}$. Recall that such a scroll may be described geometrically by fixing disjoint subspaces $\PP^{a}, \PP^{b} \subset \PP^{a+b+1}$,  rational normal curves $C_{a}\subset \PP^{a}$ and $C_{b}\subset \PP^{b}$ of degrees $a$ and $b$ respectively,
and a one-to-one correspondence $\phi\subset C_{a}\times C_{b}$.  We write $S = S_\phi$ for the correspondence scroll, which is the union of the lines $\overline{(x,y)}$ for $(x,y)\in \phi$. When $a,b\geq 1$
the surface $S$ is a smooth rational surface of degree $a+b$,   isomorphic to
$${\rm Proj}_{\PP^{1}}(\CO_{\PP^{1}}(a-b) \oplus \CO_{\PP^{1}}).$$


In addition to the double structure on $S$ that is the $K3$ carpet $X(a,b)$, we will also study the equations of a family of reducible K3 surfaces,  the union of two  scrolls
$S_1 \cup S_2$ that degenerates to $X(a,b)$. We take
$S_1 = S = S_\phi$ and
define $S_2 = S_{\phi\tau}$ as the scroll corresponding to
the correspondence $\phi\circ(\tau\times 1)\subset C_{a}\times C_{b}$, where $\tau$ is an automorphism of 
$C_a \cong \P^1$. Finally, we set 
$$
X_{\phi, \tau} = S_1\cup S_2.
$$

Now suppose that $\tau$ has two distinct fixed points, which we take to be 0 and $\infty$. In this case we may
 identify $\tau$ as multiplication by a  scalar $t\neq 1$. Had we reversed the roles of $0$  and $\infty$ (or of $C_{a}$ and $C_{b}$ we would replace $t$ by $t^{-1}$. but up to these changes $t$ is well-defined by the (abstract) surface $X_{\phi, \tau}$ as the ratio of
the points of $C_a\setminus\{0,\infty\}$ corresponding to a given point of 
$C_b\setminus\{\phi(0),\phi(\infty)\}$.
 
 The intersection 
$S_{\phi}\cup S_{\phi\tau}$ is a curve of degree $a+b+2$ and arithmetic genus 1 consisting of $C_{a}\cup C_{b}\cup L_{0}\cup L_{\infty}$, where $L_{0}, L_{\infty}$ are
 the rulings of either scroll through the points 0 and $\infty$ on $C_{a}$.
 
  
We may let $t$ go to 1, and when this happens the union of the two scrolls  approaches $X(a,b)$  (Theorem \ref{equations}).

\section{Equations and Gr\"obner Bases}\label{Sequations}
\subsection{Notation:} 
Let $a\geq b\geq 1$ be integers,  consider a projective space
$\PP^{a+b+1}_\FF$ over an arbitary field $\FF$, and let 
$$
P = \FF[x_0,x_1,\dots,x_{a},\ y_0,y_1,\dots,y_{b}]
$$
be its homogeneous coordinate ring.
Define matrices
$$
MX := 
\begin{pmatrix}
x_0&x_1&\dots&x_{a-1}\\
x_1&x_1&\dots&x_{a}
\end{pmatrix}, \qquad
MY_{t} := \begin{pmatrix}
y_0&y_1&\dots&y_{b-1}\\
ty_1&ty_2&\dots&ty_{b}
\end{pmatrix}
$$
and let
$$
M_{t} = 
\begin{pmatrix}
 x_0&x_1&\dots&x_{a-1}& y_0&y_1&\dots&y_{b-1}\\
  x_1&x_1&\dots&x_{a}& ty_1&ty_2&\dots&ty_{b}
\end{pmatrix}
$$
be their concatenation. 

We omit the subscript and write $MY$ or $M$ for $MY_{1}$ or $M_{1}$.
We will use the symbol $|$ to denote concatenation: for example, $M = MX|MY$. 

Let $I_2(MX), I_2(MY)$, and $I_2(M)$ be the ideals in $P$ generated by the $2\times 2$ minors of these matrices. In the case  $b=1$ we will also use the $2\times 2$ matrix
$$
MY2 := \begin{pmatrix}
y_0^2&y_0y_1\\
y_0y_1&y_1^2
\end{pmatrix}.
$$

Write $R := R(a,b) = P/I_2(M)$ for the homogeneous coordinate ring of the scroll $S_{t} \cong S(a,b)$ defined by $I_{2}(M_{t}).$
The line bundle corresponding to the ruling of the scroll $S_{t}$ is the cokernel of the matrix $M_{t}$, and the elements $x_0,x_1$ may be identified with the sections of this bundle. 

\subsection{The K3 Carpets}
Now let $M = M_{1} = MX|MY$.
The minimal free resolution of $I_2(M)$ is an Eagon-Northcott complex. From the form of this complex
\cite{BE} we see that
the canonical module $\omega_R$ of $R $ is isomorphic to the 
ideal 
$$
 (x_0, x_{1})^{a+b-2}R,
$$
shifted so that the generators are in degree 2, that is, 
$$
\omega_{R }\cong (x_0, x_{1})^qR(q-2).
$$

By \cite[Theorem 1.3]{GP} there exists a unique surjection
$I \to \omega_R$. We begin by making this explicit:

\begin{theorem}\label{canonical map}
Set $q = a+b-2$. The unique surjection $\alpha:I(S) \to \omega_R$ from the ideal $I(S)$ of $S$ to the module 
$ 
\omega_R
$
annihilates $I_2(MX)+I_2(MY)$ and sends
$$
\det
\begin{pmatrix}
 x_i& y_{j}\\
 x_{i+1} & y_{j+1}
\end{pmatrix}
$$
to the monomial $x_0^{q-i-j}x_1^{i+j}$.
\end{theorem}

\begin{proof} 
The given formula for $\alpha$  defines a surjection from the vector space generated by the quadrics in $I(S)$ to the vector space generated by the forms $p_\ell = x_0^{q-\ell}x_1^{\ell} \in R$.  To see that this defines a homomorphism of $P$-modules, we must show that the relations on the quadrics go to 0. 

In the case $a=b=1$ the ideal $I(S)$ is principal,  the canonical module is isomorphic to $R$, and the result is trivial. Thus we may assume that $a\geq 2$. 

The exactness of the Eagon-Northcott complex shows that the relations on the quadrics are generated by the relations on the minors of the $2\times 3$ submatrices $M'$ of $M$. Such a submatrix must involve
either two columns from $MX$ or two columns from $MY$. Since the two cases are similar, we may as well suppose that the submatrix is
$$
M' = \bordermatrix{
&0&1&2\cr
 &x_i&x_{j}&y_s\cr
& x_{i+1}&x_{j+1}&y_{s+1}
}
$$
with $0\leq i<j\leq a-1$ and $0\leq s\leq b-1$.
The relations on the minors of $M'$ are generated by
\begin{align*}
&x_i\Delta_{1,2} -x_{j}\Delta_{0,2}+y_s\Delta_{0,1} = 0\\
&x_{i+1}\Delta_{1,2} -x_{j+1}\Delta_{0,2}+y_{s+1}\Delta_{0,1} = 0.
\end{align*}
where $\Delta_{u,v}$ denotes the determinant of the $2\times 2$ submatrix of $M'$ involving the $u$-th and $v$-th columns.

The map $\alpha$ sends $\Delta_{0,1}$ to 0, so these relations go to
\begin{align*}
&-x_{j} p_{i+s} +x_{i}p_{j+s}\\
&-x_{j+1} p_{i+s} +x_{i+1}p_{j+s}.
\end{align*}
In the fraction field of $R$ we have 
$$
x_{1}/x_{0} \equiv x_{2}/x_{1}  \equiv \cdots \equiv y_{1}/y_{0} = \cdots \mod I(S).
$$
In particular, for $j = 0, \dots a$ we have
$$
x_j \equiv \bigl(\frac{x_1}{x_0}\bigr)^jx_0 \mod I(S).
$$
Thus the two binomials above are both congruent mod $I(S)$ to
$$
-\bigl(\frac{x_1}{x_0}\bigr)^jx_{0}x_{0}^{q-i-s}x_{1}^{i+s}+
\bigl(\frac{x_1}{x_0}\bigr)^ix_{0}x_{0}^{q-j-s}x_{1}^{j+s}
=0
$$
as required.
\end{proof}

\def\Pbar{{\overline P}}
\def\kbar{{\overline k}}

\subsection*{Some reducible K3 surfaces}
We now turn to the ideal of the K3 surfaces $X_{\phi, \tau}$ in the case where $\tau$ is multiplication by a scalar $t$. 
It turns out that it is convenient to  write down generators in some cases where $t$ is not defined over the ground
field $\FF$, but is the ratio $t = t_{1}/t_{2}$ of two the roots $t_{1}, t_{2} \neq 0$ of a quadratic equation $p(z) = z^2-e_1z+e_2 \in \FF[z]$. We include the possibility $\FF = \ZZ$ as well---this will be important in Section~\ref{resolutions}.  We write $e$ for the pair $(e_1,e_2)$. As we shall see, if $(e_1,e_2)\in \FF$ then the scheme $X_{\phi,\tau}$  has a model $X_{e}$ defined over $\FF$.

We  think of the $t_{i}$ as being in a fixed algebraic closure
$\overline \FF$ of $\FF$, and set $\Pbar :=\overline \FF[x_{0},\dots,x_{a},y_{0},\dots,y_{b}]$. If $t_{1}=t_{2}$, so that $t=1$ then, for simplicity, we will suppose that $t_{1} = t_{2}=1$. 

Other than the 
minors of $MX$ and $MY$, the forms that will enter into our description are defined as follows:
\begin{enumerate}
\item In the case $a,b\geq 2$ we let $J_e \subset S$ generated by the bilinear forms 
$$
 Q_{i,j}:= x_{i+2}y_{j}-e_1x_{i+1}y_{j+1}+e_2x_{i}y_{j+2} \eqno(1a) \label{bilinear1}.
$$
for $0\leq i\leq a-2$ and $0\leq j\leq b-2$. The ideal $J_e$ can be perhaps more conveniently specified as the ideal generated by the entries of the $(a-1) \times (b-1)$ matrix 
$$
 \begin{pmatrix} 
 x_0 & x_1 & x_2 \cr
 x_1 & x_2 & x_3 \cr
\vdots & \vdots & \vdots\cr
 x_{a-2}&x_{a-1}& x_a\cr
 \end{pmatrix}
 \begin{pmatrix} 
 0& 0 &  e_2\cr
 0 & -e_1 & 0 \cr
 1 & 0 & 0 \cr
 \end{pmatrix} 
  \begin{pmatrix} 
 y_0 & y_1 & \ldots &y_{b-2} \cr
 y_1 & y_2 & \ldots & y_{b-1} \cr
 y_2 & y_3 & \ldots & y_{b} \cr
 \end{pmatrix}  \eqno(1b) \label{bilinear2}.
 $$

\item In the case $a\geq 2, b=1$  we let $J_{e}$ be the ideal generated by the cubic forms
$$
Q_{i,0} := x_{i+2}y_{0}^{2}-e_1x_{i+1}y_{0}y_{1}+e_2x_{i}y_{1}^{2}
$$
for $0\leq i\leq a-2$, i.e. the entries of the $(a-1)\times 1$ matrix
$$
 \begin{pmatrix} 
 x_0 & x_1 & x_2 \cr
 x_1 & x_2 & x_3 \cr
\vdots & \vdots & \vdots\cr
 x_{a-2}&x_{a-1}& x_a\cr
 \end{pmatrix}
 \begin{pmatrix} 
 0& 0 &  e_2\cr
 0 & -e_1 & 0 \cr
 1 & 0 & 0 \cr
 \end{pmatrix} 
  \begin{pmatrix} 
 y_0^2  \cr
 y_0y_1  \cr
 y_1^2  \cr
 \end{pmatrix}. 
 $$
 
 \item Finally, in case a=b=1 we
let $J_{e}$ be the ideal generated by the quartic form
\begin{align*}
Q_{0,0} &:= x_{1}^2y_{0}^2-e_1x_0x_{1}y_0y_{1}+e_2x_0^2y_1^2 \\
&= (x_1y_0-t_1x_0y_1)(x_1y_0-t_2x_0y_1)
\end{align*}
\end{enumerate}

Set $I_e := I_2(MX)+I_2(MY)+J_e$. We will show that $I_{e}$ is the ideal of forms vanishing on $X_{\phi, \tau}$ and that $P/I_{e}$ is a Gorenstein ring with $\omega_{P/I_{e}}\cong P/I_{e}$ as graded modules, so that, in particular, $X_{e}$ is a degenerate K3 surface.

\begin{theorem}\label{equations}\label{generators} Let $\FF$ be any field.
$I_e := I_2(MX)+I_2(MY)+J_e$ is the saturated ideal of $X_e$.
\begin{enumerate}
\item If $t_{1} = t_{2} = 1$, hence $e=(2,1)$, then $I_{e}$ is the kernel of the map $\alpha$ of Theorem~\ref{canonical map},
and thus $I_{e}$ is the saturated ideal of $X_{e} = X(a,b)$.

\item Suppose that $t_{1} \neq t_{2}$. Define $2\times (a+b)$ matrices over 
$\Pbar$ by
$$
m_{\ell} := M_{t_\ell} =\begin{pmatrix}
x_0 & x_1 & \ldots & x_{a-1} &   y_0 & y_1 & \ldots &  y_{b-1} \cr
x_1 &x_2 & \ldots & x_a  & t_{\ell} y_1 &  t_{\ell}y_2 & \ldots & t_{\ell}y_b \cr
\end{pmatrix}
$$
for $\ell = 1,2$.
We have
$$
I_{e} = I_{2}(m_{1})\cap I_{2}(m_{2}) \subset \Pbar.
$$
and thus $I_{e}$ is the saturated ideal of a $\FF$-scheme $X_{e}$ that becomes isomorphic
over $\overline \FF$
 to $X_{\phi\tau}$,
 which is the union of the two scrolls defined by 
$I_{2}(m_{1})$ and $I_{2}(m_{2})$. These two scrolls meet along a reduced
curve
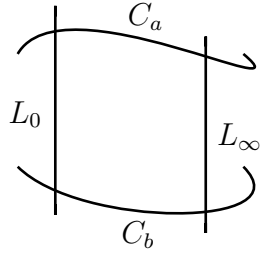
\begin{center}
\begin{tikzpicture}
\node (Ca) at (0.7,2.5) {$C_a$};
\node (Cb) at (0.6,-0.4) {$C_b$};
\node (L0) at (-0.9,1.2) {$L_0$};
\node (L1) at (1.95,0.9) {$L_\infty$};
\draw[line width=1pt] (-1,2) to [out=60, in=-40] (2,2);
\draw[line width=1pt] (-1,0.5) to [out=-45, in=-45] (2,0.5);
\draw[line width=1pt] (-0.5,0.0) to [out=-90, in=90] (-.5,2.5);
\draw[line width=1pt] (1.5,-0.25) to [out=-90, in=90] (1.5,2.1);
\end{tikzpicture} \end{center}  \noindent
where the $L_{0}, L_{\infty}$ are the lines in $\PP_{\kbar}^{a+b+1}$ defined by the vanishing of the first and
second rows of the matrix $m_{\ell}$, while the curves $C_{a}$ and $C_{b}$ are rational normal curves of degrees 
$a,b$ defined by the minors of $MX$ and $MY$ in the subspaces defined by the vanishing of the $y_{j}$ and the $x_{i}$
respectively.

\item 
The $Q_{i,j}$, together with the $2\times 2$ minors of $MX$ 
and the $2\times 2$ minors of $MY$,
form a Gr\"obner basis for $I_{e}$
with respect to the reverse lexicographic order with
 $$
 x_0>\cdots>x_a>y_0>\cdots>y_b.
 $$
\item The ring
 $P/I_{e}$ is Gorenstein, with $\omega_{P/I_{e}} \cong P/I_{e}$ as graded modules.

\end{enumerate}
\end{theorem}

We will make use of some  identities whose proofs are immediate: 
\begin{lemma}\label{identities}
Suppose that $t_{1}, t_{2}$ are nonzero scalars, and let 
$$
e_{1} = t_{1}+t_{2} \quad e_{2} = t_{1}t_{2}
$$
be the elementary symmetric functions.
\begin{enumerate}
\item If $a,b\geq 2$ then:
\begin{align*}
Q_{i,j}&:= x_{i+2}y_{j}-e_1x_{i+1}y_{j+1}+e_2x_{i}y_{j+2}\\
&= t_{2}\det
\begin{pmatrix}
 x_{i}&y_{j+1}\\
 x_{i+1} &t_{1}y_{j+2}
\end{pmatrix} 
- \det
\begin{pmatrix}
 x_{i+1}&y_j\\
 x_{i+2} &t_{1}y_{j+1}
\end{pmatrix}
\\&
= t_{1}\det
\begin{pmatrix}
 x_{i}&y_{j+1}\\
 x_{i+1} &t_{2}y_{j+2}
\end{pmatrix} 
- \det
\begin{pmatrix}
 x_{i+1}&y_j\\
 x_{i+2} &t_{2}y_{j+1}
\end{pmatrix}
\\&
\equiv
 \det
\begin{pmatrix}
x_i+t_{2}y_j& x_{i+1}+y_{j+1}\\
 t_{2}x_{i+1}+t_{1}y_{j+1} &  t_{2}x_{i+2}+y_{j+2}
\end{pmatrix}
 \mod(I_2(MX) + I_2(MY)),
\end{align*}
 
 \item If, on the other hand, $a\geq 2$ but $b=1$ then:
 \begin{align*}
Q_{i,0} &:= x_{i+2}y_{0}^{2}-e_1x_{i+1}y_{0}y_{1}+e_2x_{i}y_{1}^{2}\\
&= t_{2}\det
\begin{pmatrix}
 x_{i}&y_{0}y_{1}\\
 x_{i+1} &t_{1}y_{1}^{2}
\end{pmatrix} 
- \det
\begin{pmatrix}
 x_{i+1}&y_0^{2}\\
 x_{i+2} &t_{1}y_{0}y_{1}
\end{pmatrix}
\\&
= t_{1} \det
\begin{pmatrix}
 x_{i}&y_{0}y_{1}\\
x_{i+1} & t_{2}y_{1}^{2}
\end{pmatrix} 
-\det
\begin{pmatrix}
 x_{i+1}&y_0y_{1}\\
x_{i+2} &t_{2} y_{0}^{2}
\end{pmatrix}
\end{align*}
\end{enumerate}
 \end{lemma}\qed
 
 We will use also use the following result, which is a transposition of a well-known result on multiplicity into the context of Gr\"obner bases:
\begin{lemma}\label{AB}
Let $P = \FF[x_0\dots,x_n]$ be a standard graded polynomial ring, with a monomial order $>$, and let $I\subset P$ be a
homogeneous ideal of dimension $d$. If $g_1,\dots,g_m$ are forms in $I$ and 
$\ell_1,\dots,\ell_d$ are linear forms
such that 
$$
\length (P/(\init_<g_1,\dots,\init_<g_m, \ell_1,\dots,\ell_d)) \leq \deg P/I
$$
then $g_1,\dots, g_m$ is a Gr\"obner basis for $I$, the rings $P/I$ and $P/\init_{<}I$ are
Cohen-Macaulay,  and $\ell_1,\dots,\ell_d$ is a regular sequence modulo $\init_<I$. Moreover,
if $\sigma_{t}$, for $t\in \AA^{1}\setminus\{0\}$, is the one-parameter family of transformations of $\PP^{n}$ corresponding to the Gr\"obner degeneration 
associated to the monomial order $<$ then, for general values of $t$, the
elements  $\ell_1, \ldots, \ell_d$ form a regular sequence modulo $I_{t}.$
\end{lemma}

\begin{proof} For $t\neq 0$ we have $\deg P/\sigma_{t}I = \deg P/I$ because the transformation $\sigma_{t}$ is an automorphism of $\PP^{n}$. Moreover, by the semi-continuity of fiber dimension, $\ell_1,\dots,\ell_d$ is a system of parameters modulo $\sigma_{t}I$ for general $t$.
The degree is also semi-continuous, and $\init_{<}\sigma_{t}g_{i} = \init_{<}g_{i}$, so for general $t$, we have:
\begin{align*}
\deg P/I &= \deg P/\sigma_{t} I
\\&\leq \length P/\sigma_{t}I+(\ell_1,\dots,\ell_d)
\\&\leq
\length P/(\sigma_{t}g_{1},\dots,\sigma_{t}g_{m}, \ell_1,\dots,\ell_d)
\\&\leq  \length P/(\init_<g_1,\dots,\init_<g_m, \ell_1,\dots,\ell_d). 
\end{align*}
Our hypothesis implies that all the inequalities are equalities, 
so by \cite[Theorem 5.10]{AB} the rings $P/I$ and  $P/\init_{<}I$ are Cohen-Macaulay, and $\ell_1,\dots,\ell_d$ is a regular sequence
modulo $\init_{<}I$. 
 Since any proper
factor ring of a Cohen-Macaulay ring must have smaller degree, and since in any case $\deg \init_{<} I = \deg I$, we see that
$\init_{<} I =  (\init_<g_1,\dots,\init_<g_m)$, so $g_1,\dots,g_m$ is a Gr\"obner basis for $I$.
\end{proof}

\begin{proof}[Proof of Theorem~\ref{equations}]
It follows at once from the identities that $I_{e}$ is contained in the ideal of $X_{e}$. 

We next show that the generators of $I_{e}$ form a Gr\"obner basis.
Let $I'$ be the ideal generated by the initial forms of the generators; that is, by:
\begin{enumerate}

 \item  the initial forms of the $2\times 2$ minors of $MX$, namely
$x_{i}x_{j}$ for $1\leq i\leq j \leq a-1$;

\item the initial forms of the $2\times 2$ minors of $MY$, namely
$y_{i}y_{j}$ for $1\leq i \leq j\leq b-1$;

\item the initial forms of the $Q_{i,j}$, namely
$x_{i+2}y_{j}$ with
$0 \le i \le a-2 \hbox{ and } 0 \le j \le b-2$ if $b\geq 2$, or
 $x_{i+2}y_{0}^{2}$ with
$0 \le i \le a-2$ if $b=1$.

\end{enumerate}
Since $I'\subset \init_{<}I$, we see that $\dim S/I'\geq 3$. 
Set 
$$
P' = \FF[x_{1},\dots,x_{a}, y_{1},\dots,y_{b-1}] \cong  P/(x_0,\ x_a-y_0,\ y_b).
$$
The image of $I'$ in $P'$ contains
 every monomial of degree 2  except
$$
\{x_1y_j \mid 1\leq j\leq b-1\} \cup \{x_iy_{b-1} \mid 1\leq i\leq a\},
$$ 
every monomial of degree 3 except 
 $x_1x_a y_{b-1}$, (or $x_{1}x_{a}^{2}$ in case $b=1$),
 and every monomial of degree $\geq 4$. 
Thus $x_0, x_a-y_0, y_b$ is a system of parameters modulo $I'$ and
$P'/I'P'$ has Hilbert function $\{1,\ a+b-1,\ a+b-1,\ 1\}$. In particular,
 $$
 \dim_{k}(P'/I') = 2a+2b.
 $$ 
 
 By Lemma~\ref{AB},
 this implies that $x_0, x_a-y_0, y_b$ is a regular sequence modulo 
 $I'$ and modulo $I$; that $I' = \init_{<}I$; and that $P/I$ and $P/I'$ are Cohen-Macaulay rings of degree
 $2(a+b)$. In particular, $I_{e}$ is the saturated homogeneous ideal of $X_{e}$. This completes the proof
 of parts (1)-(3). 
 
 To complete the proof of part (4) we must show that $\omega_{P/I} \cong P/I$, and for this we
 may harmlessly assume that $\FF = \overline \FF$. In the case 
 $t_{1} = t_{2}$ this is implied by the result of
 Gallego and Purnaprajna \cite[Theorem 1.3]{GP}, so we need only treat the case $t_{1}\neq t_{2}$, where
 $X_{e} = S_1\cup S_2$ is the union of two scrolls.
 
From the fact that $P/I$ is Cohen-Macaulay, together with Hilbert function of $P/I'$, we know that the Hilbert function
of $\omega_{P/I_{e}}$ is equal to the Hilbert function of $P/I_{e}$, and it suffices to show that the annihilator of the element of
degree 0 is precisely $I_{e} = I_{2}(m_{1}) \cap I_{2}(m_{2})$. Since $\omega_{P/I_{e}}$ is a Cohen-Macaulay module, no element can have annihilator of dimension $< \dim I_{e}$; thus the annihilator of the element of degree 0 is either
$I_{e}$ or $I_{2}(m_{\ell})$ for $\ell = 1$ or $\ell = 2$. 

Now the annihilator of $I_{2}(m_{\ell})$ in $\omega_{P/I_{e}}$ is equal to $\omega_{P/I_{2}(m_{\ell})}$.
Since $S(a,b)$ is rational its canonical divisor is ineffective, so the nonzero global section of $\omega_{X_{e}}$
cannot come from either of the scrolls, and we are done.
\end{proof}

\begin{theorem}\label{rank 3}  The ideal $I(a,b)$ of the K3 carpet $X(a,b)$
contains all the  rank 3 quadrics
vanishing on the scroll $S(a,b)$, and if $a,b\geq 2$ then $I(a,b)$ is generated by them. 

The projective variety of rank 3 quadrics in $I(a,b)$ is the Veronese embedding of
$$
\nu_2: \PP\bigl(\Sym_{(a-2)}(\FF^2) \oplus \Sym_{(b-2)}(\FF^2)\bigr)
$$ 
in the subspace of 
$$
\PP\bigl(\wedge^2\Sym_{a-1}(\FF^2) \oplus \wedge^2\Sym_{a-1}(\FF^2)\bigr)
$$
spanned by the ${a+b-1\choose 2}$ rank 3 quadrics described in part (3) of Theorem~\ref{generators}.
 \end{theorem}

\begin{proof} If we identify $x_0,\dots, x_a$ with the dual basis to the monomial
basis of $\Sym_a(\FF^2)$ then we may regard $MX$ as a map from $\Sym_{a-1}(\FF^2)$ to $(\FF^2)^*$.
With this identification, writing $s,t$ for the basis of $\FF^2$, some of
the rank 3 quadrics in $I_{2}(MX)$ correspond to the $2\times 2$ submatrices of MX involving the pair of generalized columns $sf, tf$ for arbitrary
$f\in \Sym_{a-2}(\FF^2)$. We first prove by induction on $a$ that these rank 3 quadrics in $I_{2}(MX)$ generate all of $I_{2}(MX)$.
This is obvious when $a=1$. By induction we may assume that the rank 3 quadrics generate all the minors in the
first $a-1$ columns of $MX$. But for $i+1\leq a-2$ we have:
\begin{align*}
&\det\begin{pmatrix}
 x_{i} & x_{a-1}\\
 x_{i+1}& x_{a}
\end{pmatrix}\\
&=
\det 
\begin{pmatrix}
 x_{i}+x_{a-2} & x_{i+1}+x_{a-1}\\
 x_{i+1}+x_{a-1}& x_{i+2}+x_{a}
\end{pmatrix}
- 
\det\begin{pmatrix}
 x_{i} & x_{i+1}\\
 x_{i+1}& x_{i+2}
\end{pmatrix}\\
&
-
\det 
\begin{pmatrix}
 x_{a-2} & x_{a-1}\\
 x_{a-1}& x_{a}
\end{pmatrix}
+
\det 
\begin{pmatrix}
 x_{i+1}&x_{a-2}\\
x_{i+2}&x_{a-1}
\end{pmatrix}\,.
\end{align*}
All the terms on the right except the last have rank 3 and are of the given form, and the last is a minor from the first $a-1$ columns,
proving the claim.

The map from this $a+1$-dimensional space of matrices to the ${a\choose 2}$-dimensional
space of quadrics in $I_2(MX)$ is quadratic, and since the image spans $I_2(MX)$,
the map must be the quadratic Veronese embedding.

The same consideration holds for the rank 3 quadrics of $MY$. As in part (3) of
Theorem~\ref{generators},  we may obtain a further rank three quadric by adding the submatrix corresponding to
$f\in \Sym_{a-2}(\FF^2)$ to one corresponding to $g\in \Sym_{b-2}(\FF^2)$, thus giving us
a vector space $\Sym_{a-2}(\FF^2)\oplus \Sym_{b-2}(\FF^2)$ of $2\times 2$ matrices
whose determinants are rank 3 quadrics. The determinant map from this vector space to the space
of quadrics is also quadratic. Since the dimension of the space of quadrics in $I(X(a,b))$ is
${a+b-1\choose 2}$, and this space is spanned by the image of the determinant map, we see that the determinant map must be the quadratic Veronese map.

To see that $I(X(a,b))$ contains all rank 3 quadrics in $I(S(a,b))$ we do induction on $a+b$. 
If $a=b=1$, then $I(X(a,b))$ contains no quadrics, and  if $a=2,b=1$ or $a=1,b=2$ there is a unique quadric, and it does have rank 3 (Example 2 in the introduction), so the result is trivial in these cases. We now suppose that $a,b\geq 2$.

Let $Q$ be a rank 3 quadric hypersurface containing $S(a,b)$. The vertex of $Q$, which is a codimension 3 linear space, is set-theoretically the intersection of $Q$ with a general linear space of codimension 2 containing it, as one can see by diagonalizing the equation of $Q$.
Such a codimension 2 space must intersect the 2-dimensional surface $S(a,b)$, necessarily in a point $p$ lying in the vertex.
 Let $\pi: \PP^{a+b+1} \to \PP^{a+b}$
be the projection from this point. 

We may choose variables within the spaces $( x_0,\dots,x_a)$ and $( y_0,\dots,y_b)$ so that (possibly after reversing the roles
of $x,y$) the point $p$ has homogeneous coordinates
$(1,0,\dots,0)$, and thus lies on the rational normal curve $C_{a} \subset S(a,b)$. It follows that $\pi(S(a,b)) = S(a-1,b)$.

The variety $\pi(X(a,b))$ is defined by
the ideal
$$
I' := I(X(a,b)) \cap \FF[x_1,\dots,x_a,y_0,\dots,y_b],
$$
and (after renumbering the variables) this ideal  contains all the quadrics in the ideal $I(X(a-1,b))$ described in 
Theorem~\ref{generators}. Thus $\pi(X(a,b)) \subset X(a-1,b)$. Since the general codimension 2 plane through $p$ meets
$X(a,b)$ in a double point at $p$, we have $\deg\pi(X(a,b)) =  \deg(X(a,b)) - 2 = \deg X(a-1,b)$. Since 
$\pi(X(a,b))$ also has the same dimension as $X(a-1,b)$, and the latter is Cohen-Macaulay, we have $\pi(X(a,b) = X(a-1,b)$.

By induction, $X(a-1,b)$ lies on all the rank 3 quadric hypersurfaces containing $S(a-1,b)$; in particular,
it lies on $\pi(Q)$. Thus $X(a,b)$ lies on $Q$.
\end{proof}

\begin{proposition} Suppose that $t_{1}\neq t_{2}$. The scheme
 $X_{e}=S_{\phi}\cup S_{\phi\tau}$ has a transverse $A_{1}$ singularity along the intersection of the two scrolls away from the 4 double points of the curve
 $E = L_{0}\cup L_{\infty}\cup C_{a}\cup C_{b}$.
\end{proposition}

\begin{proof} We may harmlessly assume $\FF =\overline \FF$ and $a\geq b\geq 1$. Consider the affine chart $U\cong \AA^{a+b+1}$ of $\PP^{a+b+1}$ defined by $\{x_0=1\}$. This open set misses the curves
$L_{\infty}$ and $C_{b}$ that are defined by
the vanishing of the first row of the matrix $MX|MY$ and the vanishing of all the variables of $MX$, respectively.   

The variables $x_1,y_0$ restrict to global  coordinates both on
$S_{\phi} \cap U \cong \AA^2$ and $S_{\phi\tau} \cap U \cong \AA^2$. Because $0 \neq e_2\in K$, we can eliminate $x_2,\ldots,x_a$ from the coordinate ring of $X_{e}\cap \AA^{a+b+1}$ using the minors of $MX$ and, if $b\geq 2$,  we can eliminate $y_2,\ldots,y_b$ using the equations
$$
Q_{0,j}\mid_{U} = x_2y_j-e_1x_1y_{j+1}+e_2y_{i+2} \hbox{ for } j=0,\ldots,b-2.
$$
It follows that $x_1,y_0$ and $y_1$ generate the coodinate ring of the affine scheme $X_{e} \cap U$.

One remaining equation of $X_{e} \cap \AA^{a+b+1}$ in these generators is obtained from $y_1^2-y_0y_2$, which, after substitution,
corresponds
to the equation
$$
e_2y_1^2-(e_1x_1y_1-x_1^2y_0)y_0=(t_1y_1-x_1y_0)(t_2y_1-x_1y_0).
$$
All other generators reduce to zero modulo this one, since otherwise $X_{e}$ would have a component of dimension $<2$.

Thus the intersection of the two components of $X_{e}\cap U$ 
in $\AA^{3}$
defined by 
$$
y_1-\frac{1}{t_1}x_1y_0 \hbox{ and }  (\frac{t_2}{t_1}-1)x_1y_0.
$$
This set has components $x_{1} = y_{1} = 0$ corresponding to $L_{\infty}$ and $y_{0} = y_{1} = 0$ corresponding to $C_{a}$, and the
intersection is transverse away from the point $x_{0} = x_{1} = y_{1} = 0$.

The arguments for the three charts $\{x_a=1\},\{y_0=1\}\hbox{ and } \{y_b=1\}$ are similar.
\end{proof}

\section{Syzygyies over $\ZZ$ and $\ZZ/p$} \label{resolutions}

In this section we investigate the question: for which prime numbers $p$ does the
carpet $X(a,b)$ satisfy Green's conjecture over a field of characteristic $p$? We begin by unpacking this question.

Let $R$ denote a field or $\ZZ$. If $F$ is a graded free complex over a graded $R$-algebra
 with $R=P_0\cong P/P_+$ a domain, then  we set
$$
\beta_{i,j}(F):= (\rank_R F_i \tensor_P R)_{j}.
$$
Following the convention used in Macaulay2, we display the $\beta_{i,j}$  
in a \emph{Betti table}
with whose $i$-th column and $j$-th row contains the value  $\beta_{i,i+j}(F)$. If $R$ is a field or $\ZZ$ we write $X^R(a,b)$ or $X^R_e(a,b)$ to denote the the subscheme of $\PP_R^{a+b+1}$ that is defined by
the ideal described in Theorem~\ref{generators}, and we write $P^R(a,b)$ for it's homogeneous coordinate ring. 

If $F$ is the minimal free resolution  of  $P^\FF(a,b)$  as a module over
$$
\FF[x_0,\dots, x_a,y_0,\dots, y_b]
$$
where $\FF$ is a field of characteristic $p$, we say that
Green's conjecture holds for $X^\FF(a,b)$ if
$\beta_{i,i+1}(F) = 0$ for $i\geq \max(a,b)$, and similarly for $X_e^\FF(a,b)$. 
Note that the presence of the ideal of the rational normal curves of degree $a$ and $b$ inside the ideal of $X(a,b)$ implies that $\beta_{i,i+1}(F)\neq 0$ for $0< i<\max(a,b)$, so that when Green's conjecture holds, it is sharp.

We have already shown that $P^\FF(a,b)$ is Cohen-Macaulay. The hyperplane section,
which is a ribbon canonical curve, thus has minimal free resolution with the same
Betti numbers (\cite[Proposition 1.1.5]{BH93}). 
Since the hyperplane is a ribbon of genus $g=a+b+1$ and Clifford index
$b$ by \cite[p. 730]{BE} this is what Green's conjecture predicts for ribbons \cite[Corollary 7.3]{BE}.
Since ribbons do
satisfy Green's conjecture in characteristic 0 (\cite{D}), 
it follows that this is true for K3 carpets as well.
 

Returning to the general setting of 
 a graded free complex $F$ over a graded $R$-algebra $P$ with $R=P_0\cong P/P_+$, we define the $k$-th {\em constant strand} of $F$,  denoted $F^{(k)}$, to be
the submodule of elements of internal degree $k$ 
of the complex $F\tensor_P R$. Thus $F^{(k)}$ has the form:
$$ 
F^{(k)}: \cdots \leftarrow R^{\beta_{k-2,k}(F)} \leftarrow  R^{\beta_{k-1,k}(F)} \leftarrow  R^{\beta_{k,k}(F)} \leftarrow  \cdots\,.
$$ 
We write $H_i(F^{(k}))$ for the homology of this subcomplex at the term
$R^{\beta_{i,k}(F)}$. If $R$ is a field, $F$ is any graded $P$-free resolution of a module
$M$, and $F'$ is the minimal free resolution of $M$, then since the minimal free
resolution is a summand of any free resolution we have
$\beta_{i,k}(F') = H_i(F^{(k)})$.

To survey what happens for all primes $p$ at once, we  work over $\ZZ$.
We have shown that the homogenous ideal of $X(a,b)\subset \PP_\ZZ^{a+b+1}$ is minimally generated by a Gr\"obner basis consisting of forms with integer coefficients, and  the coefficients of the lead terms are $\pm 1$. Thus the homogeneous coordinate ring  $P^\ZZ(a,b)$ of $X^\ZZ(a,b) $ is a free $\ZZ$-algebra, and any free resolution over $P^\ZZ(a,b)$ reduces, modulo a prime $p$, to a free resolution of $P^{\ZZ/p}(a,b)$ over 
in characteristic $p$. 

This means that we can deduce properties in all characteristics from properties of a free resolution over $\ZZ$. We will use the (not necessarily minimal)  free resolution introduced (in a slightly different form) in \cite{S91}, called the \emph{Schreyer resolution} in Singular. See \cite{BS15} for a mathematical exposition, and \cite{EMSS} for an efficient algorithm. We have implemented a Macaulay2 package \href{https://www.math.uni-sb.de/ag/schreyer/index.php/computeralgebra}{K3Carpets.m2} \cite{ES18} for  exploration of these questions.

The definition of the Schreyer resolution of an ideal $I$,
described in  \cite{BS15},
 starts with a normalized Gr\"obner basis
$$
f_1,\ldots,f_n
$$
of $I$, sorted first by degree and then by the reverse lexicographic order of the initial terms.  
Each minimal monomial generator of the monomial ideal 
$$
M_i = ( \init(f_1),\ldots,\init(f_{i-1})) : \init(f_i) \hbox{ for } i =2, \ldots, n
$$
determines a syzygy. One shows that these syzygies form a Gr\"obner basis for the first syzygy module of $f_1,\ldots,f_n$ with respect to the induced monomial order. Their lead terms are
$m_j e_i $ for generators $e_i$ of $F_1$ mapping to $f_i$ and $m_j \in M_i$ a minimal monomial generator. 
Continuing with the algorithm, we get the finite free resolution $F$ whose terms
$F_i$ are free modules with chosen bases. 

It will be useful in the proof of Theorem~\ref{F leadTerm} to give each of the chosen basis elements of $F_p$ a name, which is a sequence $m_1,\ldots, m_p$ of monomials:

\begin{definition} The basis element $e_i$ of $F_1$ gets as a name the monomial $\init(f_i)$.
If the minimal generator $e_j \in F_p$ is mapped to a syzygy with lead term $m e_k \in F_{p-1}$, then
the name of a generator $e_j$ of $F_p$ is
$$\name(e_j) = \name(e_k),m.$$
We define the \emph{name product} of a generator $F_p$ to be the product of the monomials
in its name. The total (internal, as opposed to homological) degree of a generator is thus the degree of its name product.
\end{definition}

For simplicity, when we write $X(a,b)$, we will henceforward assume that $a\geq b$. To check whether Green's conjecture holds, we need only check a single homology group of a constant strand in an arbitrary free resolution:

\begin{proposition}\label{crucial strand} The K3 carpet $X^\FF(a,b)$ over a field $\FF$ satisfies Green's conjecture if and only if, for any graded free resolution $F$ of the homogeneous coordinate ring of $P^\ZZ(a,b)$, the constant strand $F^{(a+1)}$ satisfies $H_a(F^{(a+1)} \otimes_\ZZ \FF) =0$.
\end{proposition}
 
 \begin{proof} We must show that in the minimal free resolution $F'$ of $\P^\FF(a,b)$, the term $F'_k$, for
 $k\geq a$,  has no generators of degree $\leq k+1$. The construction of the Schreyer resolution $F$ of $P^\ZZ(a,b)$ shows that $F$ has no generators of degree $\leq k$, and since $F'$ is a summand of $F\otimes_\ZZ \FF$, the same is true for $F$. The hypothesis that that $H_a(F^{(a+1)}\otimes_\ZZ\FF)= 0$ (for any resolution $F$ over the integers) implies that $F'_a$ does not have any generators of degree $a+1$, either, proving the assertion for $k=a$. We complete the proof by induction on $k\geq a$.
 
 Assuming that 
 $F'_k$ has no generators of internal degree $\leq k+1$, the differential of $F'$ would map any generators of $F_{k+1}$ having
 internal degree $k+2$ to scalar linear combinations of generators of $F_k$ having internal degree $k+2$. Because $F'$ is minimal, this cannot happen.
%
%
\end{proof}

\begin{example}\label{X(6,6)} Here is the Betti table of the Schreyer resolution $F$ of $P^\ZZ(6,6)$ computed with Macaulay2:

\setcounter{MaxMatrixCols}{13}
\begin{small}
$$
\begin{matrix}
j \backslash i     &0&1&2&3&4&5&6&7&8&9&10&11\\ \hline
\text{0:}&1&\text{.}&\text{.}&\text{.}&\text{.}&\text{.}&\text{.}&\text{.}&\text{.}&\text{.}&\text{.}&\text{.}\\
\text{1:}&\text{.}&55&320&930&1688&2060&1728&987&368&81&8&\text{.}\\
\text{2:}&\text{.}&\text{.}&39&280&906&1736&2170&1832&1042&384&83&8\\
\text{3:}&\text{.}&\text{.}&\text{.}&1&8&28&56&70&56&28&8&1\\
\end{matrix}
$$
\end{small}



In this case, Proposition~\ref{crucial strand} shows that Green's conjecture over $\FF$ if depends only on a property of the $7$-th constant strand $F^{(a+1)}=F^{(7)}$. In our example, this has the form
$$
0\leftarrow \ZZ^8 \leftarrow \ZZ^{1736} \leftarrow \ZZ^{1728} \leftarrow 0.
$$
It has a surjective first map, so the vanishing of $ H_a(F^{(7)} \otimes_\ZZ \FF)$ is equivalent to the divisibility by $p$ of the determinant of a certain  $1728 \times 1728$ matrix $M$ over $\ZZ$. Computationally we find that
$$ \det M = 2^{1312}\, 3^{72} \, 5^{120}.$$
Thus in characteristic 0 or characteristic $p \not= 2,3,5$ this carpet satisfies Green's conjecture with Betti table
\begin{small}
$$\begin{matrix}
      &0&1&2&3&4&5&6&7&8&9&10&11\\ \hline
\text{0:}&1&\text{.}&\text{.}&\text{.}&\text{.}&\text{.}&\text{.}&\text{.}&\text{.}&\text{.}&\text{.}&\text{.}\\
\text{1:}&\text{.}&55&320&891&1408&1155&\text{.}&\text{.}&\text{.}&\text{.}&\text{.}&\text{.}\\
\text{2:}&\text{.}&\text{.}&\text{.}&\text{.}&\text{.}&\text{.}&1155&1408&891&320&55&\text{.}\\
\text{3:}&\text{.}&\text{.}&\text{.}&\text{.}&\text{.}&\text{.}&\text{.}&\text{.}&\text{.}&\text{.}&\text{.}&1\\
\end{matrix}$$
\end{small}
For the exceptional primes $p$ we can determine the Betti tables by computing the Smith normal form of $M$ and the  other matrices
in the constant strands of the non-minimal resolution.
They are\\
\begin{small}
$p=2:$
$$\begin{matrix}\hline
\text{0:}&1&\text{.}&\text{.}&\text{.}&\text{.}&\text{.}&\text{.}&\text{.}&\text{.}&\text{.}&\text{.}&\text{.}\\
\text{1:}&\text{.}&55&320&900&1488&1470&720&315&80&9&\text{.}&\text{.}\\
\text{2:}&\text{.}&\text{.}&9&80&315&720&1470&1488&900&320&55&\text{.}\\
\text{3:}&\text{.}&\text{.}&\text{.}&\text{.}&\text{.}&\text{.}&\text{.}&\text{.}&\text{.}&\text{.}&\text{.}&1\\
\end{matrix},$$
$p=3:$
$$\begin{matrix}\hline
\text{0:}&1&\text{.}&\text{.}&\text{.}&\text{.}&\text{.}&\text{.}&\text{.}&\text{.}&\text{.}&\text{.}&\text{.}\\
\text{1:}&\text{.}&55&320&891&1408&1162&48&7&\text{.}&\text{.}&\text{.}&\text{.}\\
\text{2:}&\text{.}&\text{.}&\text{.}&\text{.}&7&48&1162&1408&891&320&55&\text{.}\\
\text{3:}&\text{.}&\text{.}&\text{.}&\text{.}&\text{.}&\text{.}&\text{.}&\text{.}&\text{.}&\text{.}&\text{.}&1\\
\end{matrix}$$
\end{small}
\begin{small}
$p=5:$
$$\begin{matrix}\hline
\text{0:}&1&\text{.}&\text{.}&\text{.}&\text{.}&\text{.}&\text{.}&\text{.}&\text{.}&\text{.}&\text{.}&\text{.}\\
\text{1:}&\text{.}&55&320&891&1408&1155&120&\text{.}&\text{.}&\text{.}&\text{.}&\text{.}\\
\text{2:}&\text{.}&\text{.}&\text{.}&\text{.}&\text{.}&120&1155&1408&891&320&55&\text{.}\\
\text{3:}&\text{.}&\text{.}&\text{.}&\text{.}&\text{.}&\text{.}&\text{.}&\text{.}&\text{.}&\text{.}&\text{.}&1\\
\end{matrix}.
$$
\end{small} 
Experimentally we have strong evidence that $p=2$ and $p=5$ are also exceptional  primes for the general curve of genus $13$, while a general curve of this genus in characteristic $3$ satisfies  Green's Conjecture, see \cite{Bopp} and Remark \ref{det size} below. For characteristic  $p=2$ the experiments support the conjecture that a general smooth curve of  genus $13$ has the following Betti table with much smaller numbers
\begin{small}
$$\begin{matrix}\hline
\text{0:}&1&\text{.}&\text{.}&\text{.}&\text{.}&\text{.}&\text{.}&\text{.}&\text{.}&\text{.}&\text{.}&\text{.}\\
\text{1:}&\text{.}&55&320&891&1408&1155&64&\text{.}&\text{.}&\text{.}&\text{.}&\text{.}\\
\text{2:}&\text{.}&\text{.}&\text{.}&\text{.}&\text{.}&64&1155&1408&891&320&55&\text{.}\\
\text{3:}&\text{.}&\text{.}&\text{.}&\text{.}&\text{.}&\text{.}&\text{.}&\text{.}&\text{.}&\text{.}&\text{.}&1\\
\end{matrix},$$
\end{small}
\noindent 
then the carpet, while, for $p=5$, the experimental findings suggest that the Betti table of the carpet coincides with the conjectural Betti table of a general smooth curve of genus $13$.
\end{example}

The Schreyer resolution is rarely minimal, even for monomial ideals. Thus the following surprised us:

\begin{theorem}\label{F leadTerm} Let $a,b \ge 2$, and write $I = I_{(2,1)}$ for the saturated ideal
defining $X^\ZZ(a,b)$, as exhibited in Theorem~\ref{generators}. The Schreyer resolution of $\init(I)$  is minimal. \end{theorem}

\begin{proof} 
In our case, the minimal generators of $I$ form a Gr\"obner basis (Theorem~\ref{generators}), which is thus automatically
normalized. Let $F$ denote the Schreyer resolution of $J=\init(I)$.
Defining the $M_i$ as above, we see from the construction that the Schreyer resolution $G$ of $\init(f_1),\dots,\init(f_{n-1})$ is a subcomplex of $F$, and the quotient complex is the Schreyer resolution of $M_n$, appropriately twisted and shifted. 

There are $n={a+b-1\choose 2}$ generators of $J$, which we sort by degree refined by the reverse lexicographic order as follows
\\ $
x_1^2, x_1x_2, x_2^2, \ldots, x_{a-1}^2, x_2y_0,x_3y_0, \ldots,x_ay_0, x_2y_1,x_3y_1,\ldots,x_ay_1,y_1^2, 
x_2y_2, \\x_3y_2,\ldots,x_ay_2,y_1y_2,y_2^2,\ldots \ldots,x_2y_{b-2},x_3y_{b-2},\ldots,x_ay_{b-2},y_1y_{b-2},y_2y_{b-2},\\\ldots ,y_{b-2}^2, 
y_1y_{b-1},\ldots,y_{b-2}y_{b-1},y_{b-1}^2. 
$\\


Thus for $1 \le k \le n-1$ we have
\begin{center}
\begin{tabular}{l|c|l} 
$\init(f_k) $ & range & $M_k$ \cr \hline
$x_i x_j $ & $1 \le i \le j \le a-1$ & $( x_1,\ldots, x_{j-1} )$ \cr 
$x_i y_j$ & $2 \le i \le a-1, 0 \le j \le b-2$ & $( x_1,\ldots, x_{a-1},y_0,\ldots,y_{j-1} )$ \cr 
$x_a y_j$ & $ 0 \le j \le b-2$ & $( x_2,\ldots, x_{a-1},y_0,\ldots,y_{j-1},x_1^2 )$ \cr 
$y_i y_j$ & $1 \le i \le j \le b-2$ & $( x_2,\ldots, x_{a-1},y_1,\ldots,y_{j-1},x_1^2 )$ \cr 
$y_i y_{b-1}$ & $1 \le i < b-1$ & $( x_2,\ldots, x_{a-1},y_1,\ldots,y_{b-2},x_1^2 )$ \cr 
\end{tabular}
\end{center}
The monomial ideal $M_{n}$ is more complicated. The  initial term of $f_{n}$ is $\init(f_n)=y_{b-1}^2$, and we get
$$ M_n = ( y_1,\ldots, y_{b-2},x_1^2,x_1x_2,\ldots,x_{a-1}^2,x_2y_0,\ldots,x_ay_0 )$$ 


\begin{lemma} \label{G-res} The Schreyer resolution $G$ of the ideal   
$( \init(f_1), \ldots, \init(f_{n-1}))$ is the minimal free resolution of this ideal.
\end{lemma}

\begin{proof} For $k <n$, each $M_k$ is generated by a regular sequence of monomials.The name of each generator of $G_p$ is thus an initial monomial 
of an $f_k$, followed by an decreasing sequence of distinct elements of $M_k$ of length $p-1$.

We must show that there are no constant terms in the differential $G_{p+1}\to G_p$ for each $p>0$. The generators of $G_p$ have degrees $p+1$ and $p+2$. The $\ZZ^{a+b+2}$-grading of the monomial ideal induces a $\ZZ^{a+b+2}$-grading on $G$. Again in this grading a generator of $G_p$ has same  total degree as its name product.

Each name product of a generator of $G_p$ of degree $p+2$ is divisible by $x_1^2$ and some $y_j$. However, the only  name products of generators
of $G_{p+1}$ of degree $p+2$ that are divisible by $x_1^2$ are monomials in $\FF[x_1,\ldots,x_{a-1}]$, and the conclusion follows.
\end{proof}

To treat the case of $M_n$ we first study a smaller resolution:

\begin{lemma}\label{ell-res} The Schreyer resolution  $H$ of  the monomial ideal 
$$
J_H= ( x_1^2,x_1x_2,\ldots,x_{a-1}^2,x_2y_0,\ldots,x_ay_0 )
$$
is the minimal free resolution of this ideal.
\end{lemma}

\begin{proof} We order the monomial generators $m_k$ of $J_H$ as indicated above, and obtain
this time
\begin{center}
\begin{tabular}{l|c|l} 
$m_k $ & range & $(m_1,\ldots,m_{k-1}):m_k$ \cr \hline
$x_i x_j $ & $1 \le i \le j \le a-1$ & $( x_1,\ldots, x_{j-1} )$ \cr 
$x_i y_0$ & $2 \le i \le a-1$ & $( x_1,\ldots, x_{a-1})$ \cr 
$x_a y_0$ &  & $( x_2,\ldots, x_{a-1},x_1^2 )$ \cr 
\end{tabular}
\end{center}
As in the proof of Proposition \ref{G-res}, the generators of $H_p$ for $p\ge 1$ are in  degree $p+1$ and $p+2$, and only the name products of those in degree $p+2$ are divisible by $x_1^2y_0$, so no constant terms can occur in the differential by the $\ZZ^{a+b+2}$-grading.
\end{proof}

The resolution of $M_n$ is the tensor product of the resolution $H$ from Lemma \ref{ell-res} with the Koszul complex $\KK=\KK(y_1,\ldots,y_{b-2})$ . Thus the terms of the complex $F$ resolving $\init(I)$ are built from the terms of $G$ and terms of the tensor product complex $\KK \tensor H$ shifted and twisted:
$$
F_p = G_p \oplus \bigoplus_{q=0}^{\min(b-2,p-1)} \KK_q \tensor H_{p-1-q}(-2).
$$
Since $G$ is a subcomplex of $F$, the only possibly non-minimal parts of the differentials in $F$  have  source in the subquotient complex  $\KK(y_1,\ldots,y_{b-2}) \tensor S[-1](-2)$ and target in $G$.

The Schreyer resolution $FY$ of
$( y_1,\ldots,y_{b-1})^2$ is a  subcomplex of $F$ of which  $\KK(y_1,\ldots,y_{b-2}) \tensor S[-1](-2)$ is a subquotient. Since $FY$ has only generators of degree $p+1$ in homological degree $p \ge 1$, all maps of $FY$ and hence $F$ are minimal. This completes the proof of 
Theorem~\ref{F leadTerm}.
\end{proof}

\begin{corollary} The minimal free resolution of $\init(I)$ and the Schreyer resolution of $I$ have length $a+b-1$ and 
their non-zero Betti numbers are 
\begin{tiny}
 \begin{align*}
 \beta_{0,0}(F)\quad &=1, \cr
 \beta_{p,p+1}(F) & = p{a \choose p+1}+ \sum_{j=0}^{b-2} ((a-2) {a+j-1 \choose p-1} +{a+j-2 \choose p-1}) \cr
  & +\sum_{j=1}^{b-2} j {a+j-2 \choose p-1} +(b-2){a-2+b-1\choose p-1} \cr
  &+{b-2 \choose p-1} \quad \hbox{ for } 1 \le p \le a+b-2, \cr
     \end{align*}
      \end{tiny}
 and 
 \begin{tiny}
 \begin{align*}
 \beta_{p,p+2}(F) &=   \sum_{j=0}^{b-2} {a+j-2 \choose p-2}  
 +\sum_{j=1}^{b-2} j{a+j-2 \choose p-2} + (b-2){a-2+b-1\choose p-2}\cr 
 +  \sum_{q=0}^{p-2} {b-2 \choose q} &\big((p-q-1){a \choose p-q} +(a-p+q+1){a \choose p-q-2} + {a-2 \choose p-q-4}\big) \cr 
 &\qquad \hbox{ for }  2 \le p \le a+b-1 \cr 
  \end{align*}\end{tiny}
  \hbox{ and } 
  \begin{tiny}
  \begin{align*}
 \beta_{p,p+3}(F) & = {a+b-4 \choose p-3} \qquad \hbox{ for }  3 \le p \le a+b-1. \cr
  \end{align*}
 \end{tiny} 
\end{corollary}
 \begin{proof}
 The  complex $H$ has length $a$ and its the non-zero  Betti numbers are 
 \begin{align*}
 \, \beta_{0,0}(H)\quad &=1, \cr
 \beta_{p,p+1}(H) & = p{a \choose p+1}+(a-p) {a \choose p-1} + {a-2 \choose p-3}  \hbox{ for } 1 \le p \le a \cr
\hbox{ and } \cr 
 \beta_{p,p+2}(H) &=  {a-2 \choose p-2} \hbox{ for } 2 \le p \le a.
  \end{align*}
 The complex $G$ has length $a+b-1$ and its non-zero Betti numbers are
  \begin{align*}
 \beta_{0,0}(G)\quad &=1, \cr
 \beta_{p,p+1}(G) & = p{a \choose p+1}+ \sum_{j=0}^{b-2} ((a-2) {a+j-1 \choose p-1} +{a+j-2 \choose p-1}) \cr
  &\quad +\sum_{j=1}^{b-2} j {a-2+j \choose p-1} +(b-2){a-2+b-1\choose p-1} \cr  &\qquad \hbox{ for } 1 \le p \le a+b-2 \cr
   \end{align*}
 and 
 \begin{align*}
 \beta_{p,p+2}(G) &=   \sum_{j=0}^{b-2} {a+j-2 \choose p-2}  
 +\sum_{j=1}^{b-2} j {a+j-2 \choose p-2}  \cr & \qquad  +(b-2) {a-2+b-1\choose p-2} \hbox{ for } 2 \le p \le a+b-1. \cr
  \end{align*}
  The formula now follows  from
  $$
F_p = G_p \oplus \bigoplus_{q=0}^{\min(b-2,p-1)} \KK_q \tensor H_{p-1-q}(-2).
$$ \end{proof}

\begin{remark}\label{f(a,b)} The formula for $\beta_{p,p+1}(F)$ can be a simplified:
 \begin{align*} 
 \beta_{p,p+1}(F) &= {a-2 \choose p-1}+ {b-2 \choose p-1}
 + p{ a+b-1 \choose p+1} - 2{a+b-3 \choose p-1}.  \cr
  \end{align*}
  Using this and $
 \beta_{p-2,p+1}(F) ={a+b-4 \choose p-1} 
$
we can also obtain a simplified formula for the $\beta_{p,p+2}(F)$'s  by using the identities
  \begin{align*}
  \beta_{p,p+1}(F)&-\beta_{p-1,p+1}(F)+\beta_{p-2,p+1}(F) \cr
  &= p{a+b-3 \choose p+1}-(a+b-2-p){a+b-3 \choose a+b-1-p}\cr
  &= \frac{a+b-2-p}{p+1}{a+b-2 \choose p-1}(a+b-2p-2) .
    \end{align*}
\end{remark}

\begin{remark}\label{subcomplex} Eliminating $y_0$ from the equations of  $X_e(a,b) \subset \PP^{a+b+1}$ gives the equations of an
  $X_e(a,b-1) \subset \PP^{a+b}$, and it follows that the Schreyer resolution of $X_e(a,b-1)$ is a subcomplex of the 
 Schreyer resolution  of $X_e(a,b)$. Indeed the generators derived from
 \begin{center}
\begin{tabular}{l|c|l} 
$\init(f_k) $ & range & $M'_k$ \cr \hline
$x_i x_j $ & $1 \le i \le j \le a-1$ & $( x_1,\ldots, x_{j-1} )$ \cr 
$x_i y_j$ & $2 \le i \le a-1, 1 \le j \le b-2$ & $( x_1,\ldots, x_{a-1},y_1,\ldots,y_{j-1} )$ \cr 
$x_a y_j$ & $ 1 \le j \le b-2$ & $( x_2,\ldots, x_{a-1},y_1,\ldots,y_{j-1},x_1^2 )$ \cr 
$y_i y_j$ & $2 \le i \le j \le b-2$ & $( x_2,\ldots, x_{a-1},y_2,\ldots,y_{j-1},x_1^2 )$ \cr 
$y_i y_{b-1}$ & $2 \le i \le b-2$ & $( x_2,\ldots, x_{a-1},y_2,\ldots,y_{b-2},x_1^2 )$ \cr 
\end{tabular}
\end{center}
belong to this subcomplex. 
For the last equation with lead term $\init(f_{n'})=y_{b-1}^2$ we get  
$$ M'_{n'} = ( y_2,\ldots, y_{b-2},x_1^2,x_1x_2,\ldots,x_{a-1}^2,x_2y_1,\ldots,x_ay_1 )$$ 
which is not a subset of the corresponding $M_{n}$. Hence some generators of the Schreyer resolution for 
$X_e(a,b-1)$ are not mapped to generators
Schreyer resolution of $X_e(a,b)$ but rather to linear combinations. 
\end{remark}

\begin{remark} The equations of $X_e(a,b)$ allow a $\ZZ^3$-grading. The equations and the whole resolution is homogenous
for $ \deg x_i = (1,0,i)$ and $\deg y_j = (0,1,j)$. The non-minimal maps in the non-minimal resolution decompose into blocks 
with respect to this fine grading. 
\end{remark}

We can also compute the Betti table for the minimal resolutions of the K3 carpets $X^\FF(a,b)$ over a field $\FF$ of characteristic $2$. Note that, because $e_1,e_2$ are elements of $\FF$, the degenerate K3 surface  $X^\FF_{(0,1)}(a,b)$ coincides with the carpet
$X^\FF(a,b)=X^\FF_{(2,1)}(a,b)$.

\begin{theorem}\label{4-gonal} Let $a,b \ge 2$ and let $\FF$ be an arbitrary field .
The minimal free resolution of the homogeneous coordinate ring of  $X: = X_{e}(a,b)\subset \PP^{a+b+1}$  for $e=(0,1)$ has Betti  numbers 
$$
\beta_{i,i+1} = i {a+b-2 \choose i+1} + (\max(a-i,0)+\max(b-i,0)){a+b-2 \choose i-1}
$$ for $i \ge 1$
and $\beta_{i,i+2}=\beta_{a+b-1-i,a+b-i}$ for $1 \le i \le a+b-2$.
(These Betti numbers coincide with the Betti numbers of a $4$-gonal canonical curve of genus $g=a+b+1$ with relative canonical resolution invariants $a-2$  and $b-2$, see \cite[Example (6.2)]{S86}.
\end{theorem}

\begin{proof} The $2\times 2$ minors of the matrix
$$m=
\begin{pmatrix}
x_0 & x_1 & \ldots & x_{a-2} & y_0 & y_1 & \ldots & y_{b-2} \cr
x_2 & x_3 & \ldots & x_a      &   -y_2 &  -y_3 & \ldots & -y_b \cr
\end{pmatrix}
$$ 
are contained in $I_X$. Thus $X$ is contained in a 4-dimensional rational normal scroll of type
$$ 
Y=S(\lfloor a/2 \rfloor,\lceil a/2 \rceil-1,\lfloor b/2 \rfloor,\lceil b/2 \rceil-1)
$$
of degree $f=a-1+b-1$.
As a subscheme of the scroll, $X$ is the complete intersection of two divisors, whose classes are
of class $2H-(a-2)R$ and $2H-(b-2)R$, where $H,R \in \Pic Y$ denote the hyperplane class and the ruling of $Y$. These are defined by the vanishing of
$$
x_1^2-x_0x_2, x_2^2-x_1x_3, \ldots, x_{a-1}^2-x_{a-2}x_a
$$ 
and 
$$
 y_1^2-y_0y_2, y_2^2-y_1y_3, \ldots, y_{b-1}^2-y_{b-2}y_b, 
$$
respectively. In terms of the Cox ring $\FF[s,t,u_0,u_1,v_0,v_1]$ of $Y$
they are given by relative quadrics
$$
\begin{cases}
u_1^2-stu_0^2 \quad \hbox{ if } a \equiv 0 \mod 2 \cr
su_1^2-tu_0^2 \quad \hbox{ if } a \equiv 1 \mod 2 \cr \end{cases}
$$
and
$$
\begin{cases}
v_1^2-stv_0^2 \quad \hbox{ if } b \equiv 0 \mod 2 \cr
sv_1^2-tv_0^2 \quad \hbox{ if } b \equiv 1 \mod 2 \cr \end{cases}.
$$
Thus by \cite[Example (3.6) and (6.2)]{S86}  the minimal free resolution of $I_X$ is given by an iterated mapping cone
$$ \sC^{0} \leftarrow [\sC^{a-2}(-2) \oplus \sC^{b-2}(-2) \leftarrow \sC^{f-2}(-4) ] $$
where $\sC^j$ denotes the $j$-th Buchsbaum-Eisenbud complex associated to $m$. (The complexes
$\sC^0, \sC^1$ are also known as Eagon-Northcott complex and Buchsbaum-Rim complex of $m$.)
\end{proof}

Part of Theorem \ref{4-gonal} generalizes as follows:

\begin{theorem} \label{poncelet}{[Resonance]} Suppose $p(z)=z^2-e_1z+e_2$ has distinct non-zero roots $t_1,t_2\in \FF$ such that $t_2/t_1$ is a primitive $k$-th root of unity and $a,b \ge k+1$,
and set $X := X^\FF_e(a,b)$. 
\begin{enumerate}
\item $X$ is contained in a rational normal scroll
of type $$Y=S(a_0,\ldots,a_{k-1},b_0,\ldots,b_{k-1})$$ with
$$ 
a_i = |\{0 \le j \le a \mid j \equiv i \mod k \}| -1
$$
and
$$ b_i = |\{0 \le j \le b \mid j \equiv i \mod k \}|-1. $$
\item The map $Y \to \PP^1$ induces a fibration of $X$ into $2k$-gons.
\item If $a,b \ge 2k^2$ then $X$ has  graded Betti numbers $\beta_{\ell,\ell+1}=0$ for $\ell > a+b-1+2-2k$ and $\beta_{\ell,\ell+2}=0$ for $\ell < 2k-2$. In particular the range of non-zero Betti numbers coincides with range  predicted by Green's conjecture for a general $2k$-gonal curve of genus $g=a+b+1$.
\end{enumerate}
\end{theorem}

In characteristic 0, Green's conjecture is known to hold for general  $d$-gonal curves of every genus by \cite{A}, and it is known in every characteristic for some  $d$-gonal curve of genus $g$ if $g>(d-1)(d-2)$ by \cite{S88}. However we do not know that the family of curves of
genus $g$ and gonality $d$ is irreducible; and indeed the Hurwitz scheme could be reducible in  positive characteristics, see \cite[Example 10.3]{F69}.

\begin{proof}[Proof of Parts (1) and (2)] By Theorem \ref{equations}, $X$ is the union of the two scrolls defined by the minors of the matrices
$$
m_\ell =\begin{pmatrix}
x_0 & x_1 & \ldots & x_{a-1} &\mid&   y_0 & y_1 & \ldots &  y_{b-1} \cr
x_1 &x_2 & \ldots & x_a  &\mid& t_{\ell} y_1 &  t_{\ell}y_2 & \ldots & t_{\ell}y_b \cr
\end{pmatrix}\quad \hbox{for } \ell = 1,2$$
respectively.

Applying an automorphism of $\PP^{a+b-1}$ we may assume that $t_1=1$ and thus that
 $t=t_2$ is a $k$-the root of unity. The minors of the matrix
$$
m=\begin{pmatrix}
x_0 & x_1 & \ldots & x_{a-k} &\mid&   y_0 & y_1 & \ldots &  y_{b-k} \cr
x_k &x_{k+1} & \ldots & x_a  &\mid&  y_k &  y_{k+1} & \ldots & y_b \cr
\end{pmatrix}
$$
lie in the intersection of the ideals of minors of $m_1$ and $m_2$, as one sees from the formulas
$$\sum_{\ell=0}^{k-1} t^{k-\ell-1} 
\begin{vmatrix}
x_{i+\ell} & y_{j-\ell-1} \cr
x_{i+\ell+1} & t y_{j-\ell} \cr 
\end{vmatrix}
=
\begin{vmatrix}
x_{i} & y_{j-k} \cr
x_{i+k} & t ^ky_{j} \cr 
\end{vmatrix}\,,
$$
which hold for $0\le i \le a-k$ and $k \le j \le b$.
Thus the scheme $X$ is contained in a $2k$-dimensional scroll of the type  claimed (for example 
$$
\begin{pmatrix}
x_0 & x_k & \ldots & x_{(a_0-1)k} \cr
x_k &x_{2k} & \ldots & x_{a_0k}   \cr
\end{pmatrix}
$$
is a submatrix of $m$).

Since $X=S_1 \cup S_2$ is the union of two scrolls whose basic sections $C_a$ and $C_b$ coincide we find a pencil of $2k$-gons (away from the ramification points at 0 and infinity of the $k$-power map from $\P^1$ to $\P^1$) as follows by alternating rulings from $S_1$ and $S_2$. Starting from a general  point  $(1:s:s^2: \ldots: s^a: 0 \ldots 0) \in C_a$ we have a ruling of the first scroll $S_{1}$ connecting it to the point $(0:\ldots:0:1:s: \ldots:s^b) \in C_b$. The ruling of the second scroll $S_{2}$ joins this point on $C_b$ with the point $(1:ts:\ldots :(ts)^a:0:\ldots:0)$. 
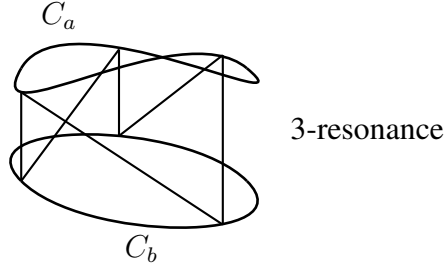
\begin{center}
\begin{tikzpicture}
\node (Ca) at (-0.5,2.7) {$C_a$};
\node (Cb) at (0.6,-0.4) {$C_b$};
\node (R) at (3.6,1.2) {3-resonance};
\draw [- ,thick] (-1,1.68) -- (-1,0.5);
\draw [- ,thick] (-1,0.5) -- (0.3,2.25);
\draw [- ,thick] (0.3,1.1) -- (0.3,2.25);
\draw [- ,thick] (0.3,1.1) -- (1.68,2.175);
\draw [- ,thick] (1.68,-0.08) -- (1.68,2.175);
\draw [- ,thick] (1.68,-0.08) -- (-1,1.68);
\draw[line width=1pt] (-1,2) to [out=60, in=-40] (2,2);
\draw[line width=1pt] (-1,2) to [out=240, in=140] (2,2);
\draw[line width=1pt] (-1,0.5) to [out=-45, in=-45] (2,0.5);
\draw[line width=1pt] (-1,0.5) to [out=135, in=135] (2,0.5);
\end{tikzpicture} \end{center}
 \noindent
Continuing with a ruling of the first scroll and so on this process closes with an $2k$-gon,
since $t$ is a primitive $k$-th root of unity. 

The map $Y \to \PP^1$ sends a point of $Y$ to the ratio of the two rows of $m$ evaluated at that point, so the $2k$-gon is contained in the fiber defined by
$$\begin{pmatrix} s^k,-1 \end{pmatrix}
\begin{pmatrix}x_0 & x_1 & \ldots & x_{a-k} &   y_0 & y_1 & \ldots &  y_{b-k} \cr
x_k &x_{k+1} & \ldots & x_a  &  y_k &  y_{k+1} & \ldots & y_b \cr
\end{pmatrix}=0,$$
Since $s^k=\tilde s$ has $k$ distinct solutions for $\tilde s \not= 0$, the fiber of the composition $X=S_{1} \cup S_{2} \hookrightarrow Y \to \PP^1$ over the point $(1:\tilde s)$ contains precisely $k$ rulings of each of the two scrolls $S_{\ell}$. Hence the $2k$-gon is the complete fiber
of $X \to \PP^1$.

The last statement follows by resolving the relative resolution of $X$ in the $2k$-dimensional scroll $Y$ by an iterated mapping cone built from Buchsbaum-Eisenbud complexes following the strategy of \cite{S88}.
Before we discuss details, we look at two examples.

\begin{example}\label{6-gonal} We consider cases of 3-resonance, $k=3$, and take $X = X_{(-1,1)}(a,b)\subset \PP^{a+b+1}_{\FF}$, since the polynomial $p(z)=z^2+z+1$ has as zeroes the primitive third roots of unity.  Note that in characteristic 3 the union of scrolls
$X_{(-1,1)}(a,b)$ coincides with the carpet $X(a,b) = X_{2,1}(a,b)$, so in characteristic 3 there is no 3-resonance, but the considerations of the free resolution below are the same.  
By \ref{generators} 
 the scheme  $X=X_{(-1,1)} \subset \PP^{a+b+1}_{\FF}$ is defined by the ideal $I_{(-1,1)}$ generated by the 
$2\times 2$ minors of the two matrices
$$
\begin{pmatrix}
x_0 & x_1 & \ldots & x_{a-1} \cr
x_1 & x_2 & \ldots & x_a \cr
\end{pmatrix}
\qquad
\begin{pmatrix}
y_0 & y_1 & \ldots & y_{b-1} \cr
y_1 & y_2 & \ldots & y_b \cr
\end{pmatrix}
$$ 
and the entries of the $(a-1) \times (b-1)$ matrix 
$$
 \begin{pmatrix} 
 x_0 & x_1 & x_2 \cr
 x_1 & x_2 & x_3 \cr
\vdots & \vdots & \vdots\cr
 x_{a-2}&x_{a-1}& x_a\cr
 \end{pmatrix}
 \begin{pmatrix} 
 0& 0 &  1\cr
 0 & -1 & 0 \cr
 1 & 0 & 0 \cr
 \end{pmatrix} 
  \begin{pmatrix} 
 y_0 & y_1 & \ldots &y_{b-2} \cr
 y_1 & y_2 & \ldots & y_{b-1} \cr
 y_2 & y_3 & \ldots & y_{b} \cr
 \end{pmatrix} 
 $$

We suppose for concreteness that $a,b \equiv 2 \mod 3$. Then the scheme $X$ is contained in a  scroll $Y$ of type 
$$Y=S(\frac{a-2}{3},\frac{a-2}{3},\frac{a-2}{3},\frac{b-2}{3},\frac{b-2}{3},\frac{b-2}{3}).$$

In terms of the Cox ring ($\equiv$ toric coordinate ring) $\FF[s,t,u_0,u_1,u_2,v_0,v_1,v_2]$ of $Y$ the remaining equations reduce to
an ideal sheaf $\sI_{Cox}$ generated by 9 relative quadrics that are the $2\times 2$ minors
of the matrices
$$
\begin{pmatrix} 
u_0 & u_1 &s u_2 \cr
 u_1 & u_2 &t u_0 \cr
\end{pmatrix}
\quad \hbox{ and } \quad
\begin{pmatrix} 
v_0 & v_1 &s v_2 \cr
 v_1 & v_2 &t v_0 \cr
\end{pmatrix}
$$
together with
$$  u_2v_0+u_1v_1+u_0v_2,\; tu_0v_0+su_2v_1+su_1v_2,\;
 tu_1v_0+tu_0v_1+su_2v_2.
$$
The relative resolution constructed in \cite[Section 3]{S86} can be regarded as a complex of free modules over the Cox ring which sheafifies to a resolution
of $\sO_X$ by locally free $\sO_Y$-modules. In our specific case it has the Betti table
\begin{small}
$$\begin{matrix}
      &0&1&2&3&4\\ 
\text{total:}&1&9&16&9&1\\ \hline
\text{0:}&1&\text{.}&\text{.}&\text{.}\\
\text{1:}&\text{.}&3&\text{.}&\text{.}&\text{.}\\
\text{2:}&\text{.}&6&16& 6&\text{.}\\
\text{3:}&\text{.}&\text{.}&\text{.}&3&\\
\text{4:}&\text{.}&\text{.}&\text{.}&\text{.}&1\\
\end{matrix}$$
\end{small} 
where we have given all the variables in the Cox ring degree $1$.

We specialise further and take $a=b=8$. Then
$$Y=S(2,2,2,2,2,2)  \subset \PP^{17}_\FF$$
is a rational normal scroll of degree $f=12$ isomorphic to $ \PP^1_\FF \times \PP^5_\FF$.

The relative resolution of $\sO_X=\sO_{X_e(8,8)}$ as an $\sO_Y$-module has shape
{\small
\begin{align*}
\sO_X \leftarrow &\sO_Y \leftarrow \sO_Y(-2H+3R)^6 \oplus \sO_Y(-2H+4R)^3 \leftarrow  \sO_Y(-3H+5R)^{16} \leftarrow \cr
 & \qquad\sO_Y(-4H+6R)^3 \oplus \sO_Y(-4H+7R)^3 \leftarrow \sO_Y(-6H+10R) \leftarrow 0
\end{align*}
}
Here $H$ and $R$ denote the hyperplane class and the ruling of $Y$.

Each term in the relative resolution is resolved by a Buchsbaum-Eisenbud complex $\sC^j$ associated to the defining matrix $m$
 of $Y$ regarded as a map $m\colon \sF\to \sG$ between vector bundles $\sF\cong \sO(-1)^f$ and $\sG \cong \sO^2$ on $\PP^{a+b+1}$.
{\small
\begin{align*}
0\leftarrow \sO_Y(jR) \leftarrow S_j\sG& \leftarrow S_{j-1}\sG\tensor \sF\leftarrow \ldots &\cr
& \ldots \leftarrow \Lambda^j \sF\leftarrow \Lambda^{j+2}\sF\tensor \Lambda^2 \sG^* \leftarrow \ldots \cr
&\qquad \ldots \leftarrow \Lambda^f \sF\tensor \Lambda^2 \sG^* \tensor (S_{f-j-2}\sG)^* \leftarrow 0 \cr,
\end{align*} }
for $0 \le j \le f-2$, see \cite{S86} and \cite[Theorem A2.10 and Exercise A2.22]{E}. Two further facts are important to us:
\begin{enumerate}
\item The complexes $\sC^j$ remain exact under the global section functor
$$\sE \mapsto \Gamma_*( \sF) = \oplus_{n\in \ZZ} H^0(\PP^{a+b+1},\sE(n)),$$
i.e. we obtain  projective resolutions of $\Gamma_*(\sO_Y(jR))$ over the polynomial ring $\FF[x_0,\ldots,x_a,y_0,\ldots,y_b] = \Gamma_*(\sO_{\PP^{a+b+1}})$. (This holds because the complexes $\sC^j$ have length $f-1 < \dim \PP^{a+b+1}$.)
\item The complex $\sC^j$ has $j$ linear maps followed by a quadratic map and further linear maps.
\end{enumerate}

By (1) we can resolve the relative resolution by the iterated mapping cone of complex $\sC^j(-d)$'s. In our specific example this is the iterated mapping cone
{\small
$$
\begin{matrix}
&&\oplus^6 \sC^3(-2)&&& \oplus^3 \sC^7(-4) && \cr
 \sC^0& \leftarrow [&  \oplus & \leftarrow[ \; \oplus^{16} \sC^5(-3)&\leftarrow [  &  \oplus & \leftarrow \sC^{10}(-6) \; ]\;]\;]\cr
&& \oplus^3 \sC^4(-2) &&& \oplus^6 \sC^6(-4) \cr
\end{matrix}
$$}
The iterated mapping cone $F$ is not minimal. However, the complex $\sC^j(-d)$ for $d\ge 2$ does not contribute to the linear strand in a range outside the  contribution of the Eagon-Northcott complex $\sC^0$, which proves assertion (3) of Theorem \ref{poncelet} in this specific case.
Indeed, the additional contribution of maximal homological degree  comes from the complex $\oplus^3 \sC^7(-4)[-3]$. It is a contribution to
 $\beta_{10,11}(F) = \dim (F_{10} \tensor_S \FF)_{11}$ to which also $\sC^0$ contributes since $$10 < \length \sC^0 =f-1=11.$$
 
The presence of $\sC^0$ and its dual inside the minimal resolution gives a lower bound on the Betti numbers, which is realized for example in the case of $X_{(-1,1)}(6,6)$ in characteristic 3 computed in Example~\ref{X(6,6)}, and therefore in characteristic 0 and all but finitely many other primes. Further computation shows  that the only exceptional primes for $X_{(-1,1)}(6,6)$ are 2 and 5.
\end{example}

\noindent {\it Proof} of Theorem \ref{poncelet} (3).
We continue with the proof of Theorem \ref{poncelet} keeping the notation of the first part of the proof.

The Cox ring $\FF[s,t,u_0,\ldots,u_k,v_0,\ldots,v_k]$ is $\ZZ^2$-graded with $s,t$ of degree $(0,1)$ and  
$\deg u_i =(1,-a_i)$ and $\deg v_i=(1,-b_i)$. 
The ideal $I_{Cox}$ of $X=X_e(a,b)$ in the Cox ring
is obtained by substituting
$$x_j = s^{a_i-\ell}t^\ell u_i \hbox{ if } j=\ell k + i \hbox{ with } 0  \le i < k$$
and 
$$y_j = s^{b_i-\ell} s^\ell v_i \hbox{ if } j=\ell k + i \hbox{ with } 0  \le i < k$$
into the generators of the ideal $I_e$ and saturating with the ideal $(s,t)$.

We can alter and refine this grading to a $\ZZ^3$-grading by setting
$\deg s=\deg t=(0,0,1)$, $\deg u_i =(1,0,a_0-a_i)$ and $\deg v_i=(0,1,b_0-b_i)$,
since the substituted equations are homogeneous with respect to this grading. The last component of the degree of each variable
of the Cox ring is now $0$ or $1$.

For the description of the generators of $I_{Cox}$ the residues $0 \le \alpha,\beta< k$ with $\alpha \equiv a, \beta\equiv b \mod k$ will play a role.
Writing $j=\ell k +i$ as above the $j$-th column of the matrix $MX$ after substitution becomes
$$
\begin{pmatrix} x_j \cr x_{j+1} \end{pmatrix}=\begin{pmatrix}  s^{a_i-\ell} t^\ell u_i \cr s^{a_{i+1}-\ell}  t^\ell u_{i+1} \end{pmatrix}
\hbox{ or } \begin{pmatrix}    s^{a_{k-1}-\ell} t^\ell u_{k-1}  \cr   s^{a_0-\ell-1}t^{\ell+1} u_0 \end{pmatrix}
$$
in case $j +1\equiv 0 \mod k$.
Thus the minors of the $2\times k$ matrix  
$$A=\begin{pmatrix} 
u_0 & u_1 & \ldots & su_\alpha & \ldots& u_{k-1} \cr 
u_1 & u_2 & \ldots & u_{\alpha+1}  & \ldots & tu_0 \cr
\end{pmatrix}$$
lie in $I_{Cox}$, where the factor $s$ occurs only once in the first row, more precisely in front of $u_\alpha$,  and the factor $t$ occurs once  in the second row in front of $u_0$. Likewise we get a $2\times k$ matrix  $B$ involving the  $v$'s.

A similar pattern arises from the $(a-1) \times 3$ and $3\times (b-1)$ Hankel matrices entering the definition of the  bilinear equations 
(\ref{bilinear2}) of $X_e(a,b)$. The Hankel matrix involving the $x$'s becomes the $(k-1)\times 3$ matrix  $A'$ which is the transpose of 
$$\begin{pmatrix} 
u_0 & u_1 & \ldots & su_{\alpha-1} & su_{\alpha} &\ldots& u_{k-2} \cr 
u_1 & u_2 & \ldots & su_\alpha &  u_{\alpha+1} &\ldots& u_{k-1} \cr
u_2 & u_3 & \ldots & u_{\alpha+1} & u_{\alpha+2} & \ldots & tu_0 \cr
\end{pmatrix}$$
There are all  together at most three factors $s$ and one factor $t$. Similarly we get a $3\times (k-1)$ matrix  $B'$ involving the  $v$'s. The generators of $\sI_{Cox}$ 
of degree $(1,1,*)$ are obtained  from the entries of the $(k-1)\times (k-1)$ matrix
$$ C=A' D B'$$
with $D$ the $3 \times 3$ anti-diagonal matrix with entries $1,-e_1,e_2$ from (\ref{bilinear2}).
The ideal generated by entries of $C$ might be not saturated with respect to $st$. For example, the form
$$su_{\alpha +1}v_{\beta-1}-e_1s^2u_{\alpha}v_{\beta}+e_2su_{\alpha -1}v_{\beta+1}$$
is divisible by $s$.

By~\cite{S86} there are exactly ${2k-1\choose 2}-1$ relative quadrics. 
From the calculation above we see ${k \choose 2}$ relative quadrics of each of types
$(2,0,*)$ and $(0,2,*)$, and $(k-1)^2$ relative quadrics of type $(1,1,*)$. Since
$$ 
2{k \choose 2} + (k-1)^2= {2k-1\choose 2}
$$
we see that there is one superfluous relative quadric, and since the ones of type
$(2,0,*)$ and $(0,2,*)$ are independent, it is of type $(1,1,*)$.
 In summary, the ideal sheaf $\sI_{Cox}$ depends only on the residue classes $\alpha,\beta$ of $a$ and $b \mod k$
and is generated by
$$ 2{k \choose 2} + (k-1)^2-1= {2k-1\choose 2}-1$$
relative quadratics of degrees $(2,0,*),(0,2,*),(1,1,*)$ where $*$ represents values between 0 and 4.

The $\ell$-th free module in our relative resolution $E_\ell$ has generators of degree $(d_1,d_2,d_3)$ with $d_1+d_2=\ell+1$ for $1 \le \ell \le 2k-3$.
The last module is cyclic with a generator of degree $(k,k,2k-\alpha-\beta)$. Indeed, this is the sum of the degree of all variables of the Cox ring, which equals the degree of the generator of its canonical module. By  adjunction the relative resolution has to end with this term, since $X_e(a,b)$  has a trivial canonical bundle. The resolution is self-dual. 

The sequences
$$
\underline d_\ell = \min \{d_3 | \exists \hbox{ a  generator of $E_\ell$ of degree } (d_1,d_2,d_3) \hbox{ with } d_1+d_2=\ell+1\}
$$
and
$$
\overline d_\ell = \max \{d_3 | \exists \hbox{ a  generator of $E_\ell$ of degree } (d_1,d_2,d_3) \hbox{ with } d_1+d_2=\ell+1\}
$$
are weakly increasing, because  for each generator of the Cox ring  the third component of its degree is non-negative.

We write $\Pic(Y )= \ZZ H \oplus \ZZ R$, where $H$ denotes a hyperplane section and $R$ a fiber of $Y \to \PP^1$. In terms of the $\Pic(Y)$-grading a generator of degree $(d_1,d_2,d_3)$ corresponds to a summand
$$\sO_Y(-(d_1+d_2)H+(d_1a_0+d_2b_0-d_3)R).$$

To establish assertion (3) of Theorem \ref{poncelet}  we must show that the multi-degree $(d_1,d_2,d_3)$ of every generator of $E_\ell$ for $1 \le \ell \le 2k-3$ satisfies
$$ d_1+d_2-1+d_1a_0+d_2b_0-d_3 \le \deg Y-1=f-1.$$
Indeed, the left hand side is the length of the contribution of 
$$\sC^{d_1a_0+d_2b_0-d_3}(-d_1-d_2)$$
 to the linear part of the iterated mapping cone,
while the right hand side is the length of the $\sC^0$. 

Note that  $-d_3 \le -\underline d_\ell =-(2k-\alpha-\beta)+\overline d_{(2k-2-\ell)}$  holds by the self-duality
of the relative resolution. Because $\omega_X \cong \sO_X$ the last term in the relative resolution has to be
$\sO_Y(-2kH+(f-2)R) \cong \omega_Y$ so $f-2= ka_0+kb_0-(2k-\alpha-\beta)$.

Thus utilizing $a_0 \ge b_0$, we see that the conditions
$$\ell - (2k-1-\ell) b_0 + \overline d_{2k-2-\ell} \le 1$$
suffice. We use the rough estimate
$\overline d_{2k-2-\ell} \le 2k$, which holds since the maximal $d_3$ in the relative resolution is $2k-\alpha-\beta \le 2k$.
The desired inequality holds for all $\ell$ with  $1 \le \ell \le 2k-3$ if
$$
 b_0 \ge  2k-2 =\max \{ \frac{2k+\ell -1}{2k-1-\ell} \mid \ell=1, \ldots, 2k-3 \}
$$
Since $b+1 =kb_0-(k-1-\beta) \le kb_0$  this follows from our assumption  $a \ge b \ge 2k^2$.
\end{proof}

\begin{remark}
A proof of Theorem \ref{poncelet} (3) for $a,b \gg k$ can be deduced by  substantially easier arguments, which do not rely on the description of $\sI_{Cox}$
but only on the existence of a relative resolution proved in \cite{S86} and an analysis of how the numerical data  change when we re-embed $Y$ by $H'=H+jR$. Since
\begin{itemize}
\item $(a,b)$ will be replaced by $(a+jk,b+jk)$ and thus $f$ by $f+2jk$ and
\item $\sO_Y(-dH+cR)=\sO_Y(-dH'+(c+dj)R)$
\end{itemize}
the  conclusion of (3) is obvious for $j$ sufficiently large.
Based on experiments we conjecture that the optimal bound  is $a\ge b  \ge k^2-k$. This is true for $k\le 5$. 
\end{remark}

For further information and conjectures about relative resolutions of canonical curves see  \cite{BH15},\cite{BH17}. 

\section{Conjectures and computational results}\label{conjectures}

\begin{remark}\label{dets} It follows from Proposition \ref{crucial strand} that Green's Conjecture is true for the balanced carpet $X(a,a)$ if and only if a certain
$f(a) \times f(a)$ integer matrix has a non-zero determinant, where $$f(a)= a {2a-1 \choose a+1}-2{2a-3 \choose a-1}$$ 
by Remark \ref{f(a,b)}. By Theorem \ref{4-gonal} we know that $\beta_{a,a+1}(X(a,a)) = a { 2a-2 \choose a+1}$ over fields of characteristic $2$. Hence
$$ 2^{a{ 2a-2 \choose a+1}}$$
is a factor of this determinant. For small $a$ the relevant values are:
\begin{center}
\begin{tabular}{|c|cccccc|}\hline
$a$ & 2 & 3 & 4 & 5 & 6 & 7 \cr 
$|\det|$ & 1 &  $2^4$ &  $2^{32}3^6$  &  $2^{266} 3^{15}$ & $2^{1312} 3^{72} 5^{120}$   & $2^{6774}3^{1020}5^{315}$\cr 
$f(a)$                           & 0 & 9      & 64                & 350                   & 1728                                       & 8085 \cr
$a{2a-2 \choose a+1}$& 0 & 3      & 24                & 140                   & 720                                         & 3465 \cr  \hline
\end{tabular}
\end{center}
One step in achieving  a proof of Green's conjecture using K3 carpets might be to give an explanation of the prime power factorizations of the determinants in the table above. 

The data in this table was produced by our  Macaulay2 \cite{M2} package \href{https://www.math.uni-sb.de/ag/schreyer/index.php/computeralgebra}{K3Carpets.m2} 
 version 0.5 \cite{ES18}. Here is, how these determinants are actually computed. The first step is the computation of the Schreyer resolution of an carpet $X(a,a)$ over $\FF[x_0,\ldots,x_a,y_0,\ldots,y_a]$ for a large finite prime field $\FF= \ZZ/(p)$. In practise we take $p=32003$. The second step is to lift  the matrices in the resolution 
 to $P=\ZZ[x_0,\ldots,x_a,y_0,\ldots,y_a]$ by using the bijection of $\ZZ/32003$ with the integers in the interval  $[-16001,16001]$.  
 The resulting matrices define the Schreyer resolution over $P$ if and only if the lifted matrices  form a complex. After checking this,  we use the fine grading to find the blocks in the crucial constant strand. For the computation of the determinants of the blocks we use  their Smith normal forms. The final step is the factorisation of the product of all determinants of all blocks. 

\end{remark}

\begin{remark}\label{det size}
The enormous size ot the determinants in Remark \ref{dets} must correspond to a combination of the resonance phenomenon with the exceptional behaviour of Green's conjecture in positive characteristic.

Experimental data of \cite{Bopp}, see also \cite{BS18}, suggests that a general canonical curve of odd genus $g=2a+1$ violates Green's conjecture
in small characteristic in the following cases:
\begin{center}
\begin{tabular}{|c|c|c|c|}\hline
$a$& $g=2a+1$ & primes & $\beta_{a-1,a+1}=\beta_{a,a+1}$ \cr \hline
3 &7 & 2 & 1 \cr
4 &9 & 3 & 6 \cr
5 & 11 &  2, 3 & 28, 10 \cr
6 & 13 & 2, 5 & 64, 120 \cr
7 & 15 & 2, 3, 5 & 299, 390, 315 \cr
\hline
\end{tabular}
\end{center}
For genus $g=7,9$ this is rigorously proven by \cite{S86} and \cite{Muk}. For genus $g=11,13,15 $ we know that 
the examples found in \cite{Bopp} violate the full Green conjecture, however we do not know whether their Betti numbers coincide
with the Betti numbers of the general curve of the given genus in these characteristics. 


Computing a non-minimal resolution of the K3 union of scrolls $X_e(a,a)$ over the coefficient ring $\ZZ[e_1,e_2]$ we find the following values of the determinant of the crucial non-minimal part
\begin{center}
\begin{tabular}{|c| c |}\hline
a & $\pm \det$ \cr\hline
3 & $2e_1^3e_2^3$\cr
4 & $ 3^6e_1^{32}e_2^{32}$\cr
5 &$ 2^{46}3^{10}e_1^{220}e_2^{235}(e_1^2-e_2)^5$ \cr
6 &$2^{64}5^{120}e_1^{1248}e_2^{1464}(e_1^2-e_2)^{72}$ \cr
7 & $2^{390}3^{390}5^{315}e_1^{6377}e_2^{8302}(e_1^2-e_2)^{630}(e_1^2-2e_2)^{7}$ \cr
\hline
\end{tabular}
\end{center}

\end{remark}

Based on these values we propose two conjectures:

\begin{conjecture} \label{char 0}
For $e=(e_1,e_2) \in \FF^2$ with $e_2\not=0$ the union of scrolls $X_e(a,a)$ has a pure resolution over an field $\FF$ of characteristic 0
unless the polynomial $p(z)=z^2-e_1z+e_2=(z-t_1)(z-t_2)$ has  roots such that $t_2/t_1\neq 1$ is a $k$-th root of unity for some $k \le \frac{a+1}{2}$.
\end{conjecture}

\begin{conjecture}\label{char p} For general $e=(e_1,e_2) \in \overline \FF^2$ the  union of scrolls  $X_e(a,a)$ over an algebraically closed field $\overline \FF$ of characteristic $p$  has a pure resolution 
if $p \ge a$. In particular, Green's conjecture holds for the general curve over a field of characteristic $p$ of  genus $g$ if $p \ge \frac{g-1}{2}$. 
\end{conjecture}

 By the table above and Remark \ref{subcomplex}  both Conjectures hold for $g \le 15$.

\bibliographystyle{ABC99}

\bigskip

%
%

\end{document}

%% file: preamble.tex
\def\antiddot{\mathinner{\mkern1mu\raise1pt\vbox{\kern7pt\hbox{.}}\mkern2mu
        \raise4pt\hbox{.}\mkern2mu\raise7pt\hbox{.}\mkern1mu}}

\newcommand{\FF}{{\mathbb F}}

\newcommand{\KK}{{\mathbb K}}
\newcommand{\PP}{{\mathbb P}}

\newcommand{\ZZ}{{\mathbb Z}}

\newcommand{\HH}{{\rm{H}}}

\newcommand{\s}{\mathcal}

\newcommand{\sC}{{\s C}}

\newcommand{\sE}{{\s E}}
\newcommand{\sF}{{\s F}}
\newcommand{\sG}{{\s G}}

\newcommand{\sI}{{\s I}}

\newcommand{\sO}{{\s O}}

\newcommand{\tensor}{\otimes}

\newcommand{\punkt}{\hspace{-.3ex}\raise.15ex\hbox to1ex{\Huge.}}

\DeclareMathOperator{\Pic}{Pic}

\DeclareMathOperator{\Sym}{Sym}

\DeclareMathOperator{\length}{length}

\DeclareMathOperator{\rank}{rank}

\newtheorem{theorem}{Theorem}[section]
\newtheorem{lemma}[theorem]{Lemma}
\newtheorem{proposition}[theorem]{Proposition}
\newtheorem{corollary}[theorem]{Corollary}
\newtheorem{conjecture}[theorem]{Conjecture}
\theoremstyle{definition}
\newtheorem{definition}[theorem]{Definition}

\newtheorem{remark}[theorem]{Remark}

\newtheorem{example}[theorem]{Example}

%% file: K3CarpetEquations.bbl
\begin{thebibliography}{ABC99}




\bibitem[A05]{A} M.~Aprodu: Remarks on syzygies of d-gonal curves Green's conjecture for curves,  Math. Res. Lett. 12 (2005), no. 2-3, 387--400.

\bibitem[AF11]{AF} M.~Aprodu and G.~Farkas: Green's conjecture for curves on arbitrary {$K3$} surfaces, Compos. Math. 147 (2011), 839--851.

 \bibitem[AFPRW]{AFPRW} M.~Aprodu, G.~Farkas, S.~Papadima, C.~Raicu and J.~Weyman: In preparation.

\bibitem[AB58]{AB} M.~Auslander and D.~Buchsbaum: Codimension and Multiplicity. Annals of Mathematics 68 (1958) 625--657.


\bibitem[BE95]{BE} D.~Bayer and D.~Eisenbud: Ribbons and their canonical embeddings. Trans. Amer. Math. Soc. 347 (1995) 719--756. 


\bibitem[BS15]{BS15}C.~Berkesch and F.-O.~Schreyer:
Syzygies, finite length modules, and random curves, in:
Commutative algebra and noncommutative algebraic geometry, {V}ol. {I},
Math. Sci. Res. Inst. Publ.,
67 (2015),
25--52.

\bibitem[B17]{Bopp} C.~Bopp: Canonical curves, Scrolls and K3 surfaces. \href{https://publikationen.sulb.uni-saarland.de/bitstream/20.500.11880/26917/1/phd_bopp.pdf}{Dissertation, Saarbr\"ucken Fall 2017}.

\bibitem[BH15]{BH15} C.~Bopp and M.~Hoff: Resolutions of general canonical curves on rational normal scrolls. Archiv der Mathematik (Basel) 105 (2015) 239--249.

\bibitem[BH17]{BH17} C.~Bopp and M.~Hoff: Moduli of lattice polarized K3 surfaces via relative canonical resolutions. Preprint \href{https://arxiv.org/abs/1704.02753}{arXiv:1704.02753}.




\bibitem[BS18]{BS18} C.~Bopp and F.-O.~Schreyer:  A version of Green's conjecture over fields of finite characteristic. Preprint \href{https://arxiv.org/abs/1803.10481}{arXiv:1803.10481}.

\bibitem[BH93]{BH93} W.~Bruns and J.~Herzog, J\"urgen:
Cohen-{M}acaulay rings
Cambridge Studies in Advanced Mathematics 39, Cambridge University Press. xi, 403 p. (1993). 

\bibitem[D15]{D} A.~Deopurkar: The canonical syzygy conjecture for ribbons. \href{https://arxiv.org/abs/1510.07755}{arXiv:1510.07755}

\bibitem[DFS16]{DFS} A.~Deopurkar, M.~ Fedorchuk and D.~Swinarski: Toward GIT stability of syzygies of canonical curves.
 Algebr. Geom.  3  (2016),  no. 1, 1--22.

\bibitem[E97]{E} D.~Eisenbud: Commutative Algebra with a View toward Algebraic Geometry. GTM 150, Springer-Verlag NY, 1997.

\bibitem[E05]{E05} D.~Eisenbud: The Geometry of Syzygies. GTM 229, Springer-Verlag NY, 2005.

\bibitem[EH87]{EH} D.~Eisenbud and J.~Harris: Varieties of Minimal Degree (a centennial account). Algebraic geometry, 
Proc. Sympos. Pure Math., 46, Part 1, pp. 3-13. Amer. Math. Soc., Providence, RI, 1987. 



\bibitem[EiSa]{EiSa} D.~Eisenbud and A.~Sammartano: Correspondence Scrolls. In preparation.

\bibitem[ES18]{ES18} D.~Eisenbud and F.-O.~Schreyer: K3Carpets , a Macaulay2 package to investigate K3 carpets.
Available at\\ \href{https://www.math.uni-sb.de/ag/schreyer/index.php/computeralgebra}{https://www.math.uni-sb.de/ag/schreyer/index.php/computeralgebra}.

\bibitem[EMSS16]{EMSS} B.~Erocal, O.~Motsak, F.-O.~Schreyer and A.~Steenpass: 
Refined Algorithms to Compute Syzygies. J. Symb. Comput. 74 (2016), 308--327.

\bibitem[F69]{F69} W.~Fulton: Hurwitz schemes and irreducibility of moduli of algebraic curves. Ann. Math. (2)
90 (1969), 542--575.

\bibitem[GP97]{GP}  F.J.~Gallego and B.~Purnaprajna, Degenerations of K3 surfaces in projective space.
Trans. Amer. Math. Soc. 349 (1997) 2477--2492.

\bibitem[M2]{M2} D.R.~ Grayson and M.E.~ Stillman,
          Macaulay2, a software system for research in algebraic geometry.
          Available at \href{https://faculty.math.illinois.edu/Macaulay2/}{https://faculty.math.illinois.edu/Macaulay2/}
      
\bibitem[M95]{Muk}   S.~Mukai, Curves and symmetric spaces I. Amer. J. Math 117 (1995), 1627--1644.

\bibitem[S86]{S86} F.-O.~Schreyer: Syzygies of canonical curves and special linear series. Math. Ann. 275 (1986), 105--137.

\bibitem[S88]{S88} F.-O.~Schreyer: Green's conjecture for the general p-gonal curve of large genus. In: Algebraic curves and projective geometry, Proceedings Trento 1988, Springer Lecture Notes 1389, 254--260.

\bibitem[S91]{S91} F.-O.~Schreyer: A standard basis approach to syzygies of canonical curves. J. reine angew. Math. 421 (1991), 83--123.


\bibitem[V05]{V05} C.~Voisin: 
Green's canonical syzygy conjecture for generic curves of odd genus. Compos. Math. 141 (2005) 1163--1190. 

\end{thebibliography}
